\documentclass[final,3p,times]{elsarticle}
\makeatletter
\def\ps@pprintTitle{%
 \let\@oddhead\@empty
 \let\@evenhead\@empty
 \def\@oddfoot{}%
 \let\@evenfoot\@oddfoot}
\makeatother

\usepackage{amssymb}
\usepackage{amsthm}

\makeatletter
\newcommand*{\rom}[1]{\expandafter\@slowromancap\romannumeral #1@}
\makeatother

\def\mud{{\mu_{\rm damage}}}
\def\mun{{\mu_{\rm nuisance}}}
\def\muc{{\mu_{\rm control}}}
\newcommand{\calP}{\ensuremath{\mathcal{P}}}
\def\calPd{\calP_{\rm damage}}
\def\calPn{\calP_{\rm nuisance}}
\def\calPc{\calP_{\rm control}}
\newcommand{\RR}{\ensuremath{\mathbb{R}}}
\def\rhod{\rho_{\mud}}
\def\rhon{\rho_{\mun}}
\def\rhoc{\rho_{\muc}}

\newcommand{\calF}{\ensuremath{\mathcal{F}}}

\def\nfea{{n_{\rm features}}}
\newcommand{\EE}{\ensuremath{\mathbb{E}}}
\newcommand{\NN}{\ensuremath{\mathbb{N}}}
\def\fphy{f^{\rm exp}}
\def\fphy{\calF^{\rm exp}}
\newcommand{\calN}{\ensuremath{\mathcal{N}}}
\def\Xitt{\Xi^\text{t-t}}

\def\ntt{n^\text{t-t}}
\def\mutt{\mu^\text{t-t}}
\newcommand{\PP}{\ensuremath{\mathbb{P}}}

\newcommand{\calM}{\ensuremath{\mathcal{M}}}
\newcommand{\calC}{\ensuremath{\mathcal{C}}}
\newcommand{\calT}{\ensuremath{\mathcal{T}}}

\newcommand{\calA}{\ensuremath{\mathcal{A}}}

\newcommand{\calQ}{\ensuremath{\mathcal{Q}}}

\newcommand{\Mtens}{\underline{\calM_h}}

\newcommand{\Atens}{\underline{\calA_h}}
\newcommand{\Ahattens}{\underline{\hat\calA_h}}
\newcommand{\Ctens}{\underline{\calC_h}}
\newcommand{\Ttens}{\underline{\calT_h}}

\newcommand{\Qtens}{\underline{\calQ_h}}

\newcommand{\MRB}{\underline{\calM_{\rm RB}}}
\newcommand{\ARB}{\underline{\calA_{\rm RB}}}
\newcommand{\CRB}{\underline{\calC_{\rm RB}}}
\newcommand{\TRB}{\underline{\calT_{\rm RB}}}
\newcommand{\XRB}{\underline{Z_{\rm RB}}}
\newcommand{\QRB}{\underline{\calQ_{\rm RB}}}

\newcommand{\UU}{\ensuremath{\mathbb{U}}}

\newcommand{\bD}{\ensuremath{\mathbf{D}}}

\def\ii{\textrm{i}}
\def\ntt{n^\textrm{t-t}}
\def\ufreq{\hat{u}}
\def\mufreq{\tilde{\mu}}
\def\Pfreq{\tilde{\calP}}
\def\ufe{u_h}
\def\ufef{\hat{u}_h}
\def\ufefd{u_{h,\Delta t}}
\def\dufefd{\dot{u}_{h,\Delta t}}
\def\ddufefd{\ddot{u}_{h,\Delta t}}
\def\upr{\hat\UU_{h,\bD}}
\def\upri{\hat U_{h,\bD,i}}
\def\upro{\hat U_{h,\bD,1}}

\def\urb{U_{h,\bD,N}}
\def\urbx{U_{h,\bD,N;1}}
\def\urby{U_{h,\bD,N;2}}
\def\urbk{U_{h,\bD,N;k}}
\def\urbfd{U_{h,\bD,N,\Delta t}}
\def\urbfderr{U_{h,\bD,N,2\Delta t}}
\def\durbfd{\dot{U}_{h,\bD,N,\Delta t}}
\def\ddurbfd{\ddot{U}_{h,\bD,N,\Delta t}}

\def\afreq{\hat{a}}
\def\ffreq{\hat{f}}

\usepackage{amsmath,amsfonts,stmaryrd,latexsym,amstext}
\usepackage{bbm}
\usepackage{float}
\usepackage{algorithm}
\usepackage[noend]{algpseudocode}
\usepackage{xcolor}
\journal{Computers \& Structures}

\begin{document}

%\linespread{1.5}
%\openup 1em

\begin{frontmatter}



% \title{A Model-Order-Reduction and Neural Network-Based Approach for Damage Detection for Large Deployed Structures with Localized Operational Excitations\tnoteref{t1}}
 \title{Model-Order-Reduction Approach for Structural Health Monitoring of Large Deployed Structures with Localized Operational Excitations\tnoteref{t1}} 
 \tnotetext[t1]{This work was supported by the Office of Naval Research [N00014-17-1-2077]; and the Army Research Office [W911NF1910098]}


\author[1]{Mohamed Aziz BHOURI\corref{cor1}\fnref{fn1}%
  }
\ead{bhouri@mit.edu}

 \cortext[cor1]{Corresponding author}
 \fntext[fn1]{Present address: 3401 Walnut St, Wing A, Office 536, Philadelphia, PA 19104, USA  \\ \\ \\ Preprint submitted to Computers \& Structures}
 \address[1]{Department of Mechanical Engineering, Massachusetts Institute of Technology, 77 Massachusetts Avenue, Cambridge, MA 02139 USA}
  
\begin{abstract}
We present a simulation-based classification approach for large deployed structures with localized operational excitations. The method extends the two-level Port-Reduced Reduced-Basis Component (PR-RBC) technique to provide faster solution estimation to the hyperbolic partial differential equation of time-domain elastodynamics with a moving load. Time-domain correlation function-based features are built in order to train classifiers such as artificial neural networks and perform damage detection. The method is tested on a bridge example with a moving vehicle (playing the role of a digital twin) in order to detect cracks' existence. Such problem has $45$ parameters and shows the merits of the two-level PR-RBC approach and of the correlation function-based features in the context of operational excitations, other nuisance parameters and added noise. The quality of the classification task is enhanced by the sufficiently large synthetic training dataset and the accuracy of the numerical solutions, reaching test classification errors below $0.1\%$ for disjoint training set of size $7\times10^3$ and test set of size $3\times10^3$.
\end{abstract}



\begin{keyword}
Structural health monitoring \sep simulation-based classification \sep model order reduction \sep domain decomposition \sep parametrized partial differential equations \sep neural networks 


\end{keyword}

\end{frontmatter}


\section{Introduction}
\label{sec1}

%``model-based"
Structural Health Monitoring (SHM) has received a great deal of attention in the civil engineering and machine learning communities. The ultimate goal of SHM is to automatically identify damage before failure occurs for a given system. One alternative to solve this task is the model-based approach \cite{Khatir2019}, in which an inverse problem is solved to determine all parameters of a model from which the state of damage is inferred. The data-based approach is another predominant method in current research and consists of two stages. First, in an offline stage of a classification task, and before the structure is put into service, a set of training data is collected from a rich representation of possible healthy and unhealthy states of interest. This training dataset can be constructed (i) by performing physical experiments \cite{Farrar2013}, or (ii) by performing synthetic experiments using a mathematical best-knowledge model \cite{Hurtado2004,Hurtado2003,Lecerf2015}. Machine learning algorithms are then applied to the training dataset to obtain a classifier which assigns data to the relevant diagnostic class label. Second, in the online stage of the classification task, and during the normal operation of the structure, the classifier is used to monitor the structure, and to map measured data in the field to a best prediction of the corresponding state of damage. Simulation-Based Classification (SBC) refers to the particular choices of a data-based approach and use of synthetic experiments.

SBC for SHM of large deployed mechanical structures, such as offshore platforms and bridges, is still an active research field for several challenges associated with such task. First, such systems are subject to operational excitations \cite{Peeters2001,Moughty2017} and the most difficult loads to simulate are those with the shortest time span. By consequence, harmonic analyses fail to faithfully capture the response of these mechanical structures and time-domain characterization is needed. Indeed, modal analyses are typically not very good for local inhomogeneities since eigenfunctions do not well represent local forces whose frequency spectra are broadly spread over certain intervals. Second, numerical simulations of large structures are often challenging due to the considerable amount of memory and computation capacity needed, which is more hampering for systems with localized excitations. Such problems involve different scales, and thus numerical methods need to be accurate enough to well approximate the source terms, but also not be computationally prohibitive, which is challenging due to the large geometric domain of the structure. In this context, these forces are applied on regions which are relatively small compared to the size of the global domain. Finally, the response of the structure is inherently affected by probabilistic nuisance parameters that should be considered in the training task. Nuisance parameter refers to any parameter that affects the behavior of the structure but does not influence its state of damage. As a consequence, numerical models should consider the inherent uncertainty in the value of the parameters representing material properties, geometric domains and input forces. These deviations may be introduced in the manufacturing process or caused subsequently by operational and environmental conditions. The list of possible variations | present perhaps in both undamaged and damaged states | is thus relatively long. Moreover, ambient operational excitation will certainly deviate from ideal in terms of spatial structure and temporal signature. Therefore, it is of interest to consider not just active systems with probabilistic nuisance parameters in which one can provide forced and controlled input excitations, but also passive systems which rely on ambient loading as naturally arises in operation of the deployed system \cite{Farrar2013,Deraemaeker2008}. The goal of this work is to develop a SBC approach for such problems.

The SHM literature proposes a variety of output-only methods for passive systems. In particular, Operational Modal Analysis (OMA) identifies modal properties of structures from ambient vibration data \cite{Au2013, Moughty2017, Zhang2005,Gillich2019}. Numerous modal-based damage detection techniques have been developed. These approaches can rely on: structure's natural frequencies; modal damping; modal shape of displacement or its curvature; modal strain energy or modal flexibility. Representative examples of OMA techniques include peak-picking \cite{Maia1997}, Frequency Domain Decomposition \cite{Brincker2000,Jiang2007,Sanchez2015}, and time-domain Decomposition \cite{Karbhari2009}. All OMA techniques have been shown to be considerably sensitive to noise contamination \cite{Cruz2008,Talebinejad2011,Fan2011}. Moreover, the assumption of stationary random signals for excitations (typically white noise) made by most frequency domain techniques may not always be suited to real world applications. For instance, vehicle-induced excitations on damaged bridges can be highly non-stationary. Furthermore, the mode shapes cannot be identified precisely for complex structures. Although techniques based on modal shape of displacement or its curvature showed superior damage sensitivity compared to other OMA methods, such techniques require many sensors to recover higher modes and their performance heavily depends on the number of modes considered \cite{Farrar2013}. In addition, computing curvatures from vibration data inherently introduces additional errors due to the application of the finite difference approximation method, which is further amplified for high-frequency noise. Existing frequency domain-based classification approaches consider synthetic datasets with size that goes up to the order of $10^4$ sample points \cite{Taddei2018} but as stated earlier, frequency analyses fail to capture the structure's response to operational and localized excitations. One of the alternatives or improvement to OMA techniques is an efficient utilization of machine learning algorithms for damage detection. Recently, different vision-based methods for detecting concrete cracks have been developed \cite{Mei2019,Zhang2019,Nayyeri2019,Ni2019}. Most of these techniques rely on deep architectures of neural networks such as convolutional networks \cite{Cha2018,Li2019} and encoder-decoder networks \cite{Bang2019}. The SBC approach presented in this work is an alternative to efficiently incorporate machine learning techniques into SHM for large mechanical structures under localized operational excitations \cite{Rafiei2017}. % without calculating the defect features

Appropriate choice of features is absolutely crucial for classification. Features which are sensitive to the anticipated damage but relatively insensitive to nuisance variables and measurement noise greatly simplify the classification task and ultimately improve the robustness and hence performance of the deployed classifier. Within those considerations, structural damage detection methods using two-point time-domain correlation functions of vibration response under stochastic excitation have been developed \cite{Zhang2014,Yang2009,Huo2016}. The integrand of the correlation function clearly depends on the wave speed since it counts for two displacements with a shift in time. Therefore, it is expected to be sensitive to anticipated damage such as loss of stiffness or crack existence, since these instances of damage considerably affect the wave speed. Multiple damage detection methods have been developed based on the correlation functions such as the correlation function amplitude method (CCFA) \cite{Huo2016} | also named correlation function amplitude vector method (CorV) \cite{Yang2007} | and the inner product vector method (IPV) one \cite{Yang2009}.  The effectiveness of these methods has been shown for white noise and steady random excitations with specific frequency spectrum by visual inspection of the damage index. The latter is defined as the difference between the considered vector (CCFA or IPV) of the intact and damaged structures. Moreover, these methods were tested on relatively simple structures: the CCFA method \cite{Huo2016} was applied to a structure with $12$ degrees of freedom, while the IPV one was applied in to a clamped beam with $8$ degrees of freedom \cite{Yang2007} and to a cantilever beam with $25$ eight-node quadrilateral shell elements \cite{Yang2009}. These resolutions were justified by the necessity of obtaining rapid calculation of the vibration responses. In addition, these techniques were only tested on structures with loss of stiffness as an instance of damage. No topological changes, such as the existence of crack, were considered to further investigate their efficiency. Finally, these methods rely on the choice of a reference sensor so as to define the corresponding correlation-based vector to monitor in order to detect the damage. Similar damage detection approaches based on correlation functions have been built. The AMV and CZV methods \cite{Zhang2014} were developed as an attempt to solve the reference sensor problem. However, similarly to the CCFA and IPV methods, these techniques were only tested on relatively simple problems (for instance a structure with only $8$ degrees of freedom), did not account for nuisance parameters and also rely on visual inspection of the damage index without generalization to a SBC task.

One of the main challenges associated with the application of SBC is the construction of a sufficiently large dataset required to train the classifier. Indeed, a prerequisite for good classifier performance is a sufficiently rich description of undamaged and damaged states. Absent such a complete description (a) a classifier will certainly not be able to discriminate between different states of damage (not represented in the training dataset) and (b) the method may not be able to identify features which can discriminate between undamaged and damaged states. In particular, relying on a classical Finite Element (FE) approach to construct the dataset will be prohibitively costly. Parametric Model Order reduction (pMOR) is a mathematical and computational field of study which aims to reduce the computational cost associated with the estimation of the solution to a parameterized mathematical model. Within our context, the latter consists in the time-domain elastodynamics partial differential equation (PDE). pMOR approaches intrinsically account for the probabilistic behavior of the parameters governing the system. Therefore, the corresponding numerical model considers the inherent uncertainty in the value of the parameters. Existing pMOR techniques for time-dependent problems include Proper Orthogonal Decomposition (POD) approaches \cite{KunischVolkwein,KunischVolkwein2,RathinamPetzold}, Greedy methods \cite{Grepletal,GreplPatera}, hybrid approaches combining POD (in time) and Greedy procedures (in parameter space) \cite{Grepl}, and space-time approaches \cite{MasaThesis}.

The dimensionality of the parameter domain which describes any sufficiently rich set of possible systems | note that the parametrization characterizes material properties, geometry, boundary conditions, and also topology | precludes application of the classical pMOR approaches due to the well-known curse of dimensionality. Within the context of time-domain elastodynamics of large geometric domains with localized excitations, the recently developed two-level PR-RBC method \cite{Bhouri2020} is of great interest. The latter takes advantage of domain decomposition techniques \cite{HuynhetalESAIM,HuynhetalCMAME,EftangPatera,Smetana,SmetanaPatera} and the frequency-time duality \cite{Antoulasetal,Beattieetal} to construct the reduced bases. Using domain decomposition techniques define components which correspond to subdomains that form the global domain when they are assembled. Hence, the two-level PR-RBC method solves the curse of dimensionality issue by reducing the effective dimensionality of the parameter spaces considered in the variational problems to approximate. Therefore, it provides sufficiently accurate approximations with relatively low dimensional approximation spaces, which considerably reduces the computational cost of constructing the training dataset in the context of SBC. Moreover, the two-level PR-RBC approach addresses not a particular global system but a family of systems that can be built using a predefined library of archetype components. Such global systems are formed by connecting compatible replica of the archetype components. Hence, it offers great flexibility in topology and geometry. As a consequence, the method is particularly well-suited to SBC in which the dimension of the model parametrization is perforce very large, and where topological damages need to me modeled. 

Our goal here is to develop an efficient simulation-based classification approach for large deployed structures with localized operational excitations by considering a numerical example of a bridge with a moving vehicle. For such example, the parameters are all chosen based on literature review of actual deployed bridges. In order to overcome the previously listed challenges associated with such task, we extend the two-level PR-RBC method \cite{Bhouri2020} to account for moving loads. Moreover, using the two-point correlation functions, we build time-domain-based features which are not only sensitive to the damage | already guaranteed by the definition of the correlation function | but also relatively insensitive to nuisance parameters and measurement noise. These properties are obtained by using normalization and exploring different strategies for the choice of the time shift; the sensors; and the direction of the displacement considered for the correlation function. As for the classifiers, distinct state-of-the-art learning algorithms were tested including Artificial Neural Network (ANN) and one-vs-all Support Vector Machine (ova-SVM). The bridge example will be considered to demonstrate (a) the characterization of localized operational excitation in terms of nuisance parameters, (b) the merits of time-domain-based correlation function features in the context of localized operational excitation, other nuisance parameters and added noise for damage localization and (c) the importance of the two-level PR-RBC approach. The latter allows flexibility in topology and geometry, permits the construction of dataset with sufficiently large size and reasonable computation time, and guarantees the accuracy of the numerical solution. These last two criteria need to be satisfied in order to obtain good classification results. The methodology presented in this work can also be applied to a variety of other deployed mechanical systems such as airframes, railroad tracks, ships, \ldots Our specific contributions can be summarized in the following points:

\begin{enumerate} \item The PR-RBC-based classification permits the removal of reference sensor choice dilemma since one can actually consider all the possible combinations of correlation functions that can be built from a network of sensors and imbed them into the feature that will be used for the classification. Classifiers take care of learning from these features even if they are high dimensional.

\item Instead of visual inspection of the damage index, the PR-RBC-based method appeals to classifiers to discriminate between the different damage cases and thus provides an automated method for damage detection without human inspection. Moreover, most of the existing damage identification methods based on correlation function rely on the difference between the damage index of the undamaged and damaged state of the structure, which makes their deployment on real word systems a hard task. Supervised classification solves this issue since for any state of the structure, a feature that only depends on its current state is considered, and no reference state is needed. 
%Due to the probabilistic nature of the parameters governing a given structure's behavior, it is impossible to define a damage index corresponding to its undamaged state. Such damage index is needed as a reference within the existing damage identification methods based on correlation function cited above. 

\item The PR-RBC-based classification relies on the model order reduction technique to construct a sufficiently large dataset with a reasonable computational cost and a rich description of undamaged and damaged states by reducing the degrees of freedom from the order of $10^3-10^4$ to the order of $10-10^2$. This allows good classification performance since a classifier will be able to discriminate between different states of damage (represented in the rich training dataset) and to identify features which can discriminate between undamaged and damaged states.

\item The PR-RBC-based approach addresses full elasticity, which is carried out sufficiently fast thanks to the use of model order reduction, and in particular the two-level PR-RBC method. The method can be applied for 2-d and 3-d elasticity problems. For the current work, only the 2-d case is considered. As a consequence, the system is modeled as a continuum and the elastodynamics PDE governing its behavior is solved. This provides a more faithful characterization of the structure than modeling it as an assembly of rigid beams as it is the case for most of SHM systems detailed in the literature. The rigid beam model fails to capture effects related to cracks and other important ``health" issues such as local parameters variation but was chosen in most of the previous works due to the cheap computational cost related to solving the corresponding mathematical model with few degrees of freedom. Such local topological damages and other local parameters' variation can be considered in the elasticity model considered in this work, and the two-level PR-RBC approach inherently permits to account for this probabilistic nature. %The Model Order Reduction approach considerably reduces the degrees of freedom from the order of $10^3-10^4$ to the order of $10-10^2$, and hence significantly reduces the computational cost of the construction of the sufficiently large dataset with a sufficient accuracy to guarantee good classification results.

\item The PR-RBC-based classification enables including nuisance parameters which are inherently present in actual deployed systems such as the temporal and spatial signature of loads, damping coefficients and materials properties which are all stochastic quantities and only determined within a certain interval of confidence. The references cited above do not consider nuisance parameters since they rely on visual inspection of the damage index and not on SBC approach. The model order reduction technique and the classification task inherently permit to account for the probabilistic nature of the nuisance parameters.

\item Due to its great flexibility in topology and geometry, the PR-RBC method can accommodate for more realistic damage instances such as crack existence and is thus particularly well-suited to SHM. In such context, the computational framework should be able to switch in/out damage instances efficiently, including topological changes like cracks, and should be able to treat local variations including local excitations with support very small compared to the size of the structure dimension which is the case for bridges in instance with moving vehicles. The PR-RBC approach also offers flexibility with respect to re-use of components either for variations on a given structure or even for new but similar structures built from the same library of archetype components considered for the two-level PR-RBC method.

\item Finally, relying on the correlation function-based features, the bridge example considered in this work shows that only a sparse configuration of sensors can be required to reach sufficiently small classification errors for damage existence identification.

\end{enumerate}
%,and we construct our features inspired by those methods but adapted to the problems we consider. 

%\begin{figure} [H]
%	\begin{center}
%	\includegraphics[scale=0.11]{bridge_model_w_real_bridge_new_comp}
%	\caption{Deployed bridge and a replicating model considered for SBC for SHM}\label{fig:bridge_model_w_real_bridge_new_comp}
%	\end{center}	
%\end{figure}

%In addition to correlation functions-based method, the SHM literature proposes a variety of approaches for passive systems | ``output-only" methods | which rely on ambient loading as naturally arises in operation of the deployed system: in particular we can cite Operational Modal Analysis (OMA). Among the numerous modal-based damage detection techniques, we can cite: the natural frequencies-based approaches (Adams et al., 1978) , (Pau et al., 2011), where frequency shifts are used to detect damage; modal damping-based techniques which relies on the increase of internal friction due to cracking, resulting in a raise of the section's damping (Yamaguchi et al., 2013);  and the modal shape (Zhao \& Zhang, 2012) , (Balsamo et al. 2013) and modal curvature-based methods (Abdel Wahab \& De Roeck, 1999), which rely on node displacement and their second derivative in space respectively to detect and locate damage.
 
This paper is organized as follows. In Section \ref{sec:2levelPRRBC}, we present an overview of the recently developed two-level PR-RBC method and how to extend it in order to apply it for elastodyanmics PDE with moving loads. In Section 
\ref{sec2:SBC}, we define the Simulation-Based Classification task and the correlation function-based features considered in this work. In order to demonstrate the capability and assess the performance of the proposed technique, a numerical example is presented in Section \ref{sec:num_ex}. Finally, in section \ref{sec:Cl} we summarize our key results, discuss the limitations of the proposed approach, and carve out potential directions for future investigation.

\section{Two-Level PR-RBC Method}\label{sec:2levelPRRBC}

In this section we present a summary of the two-level PR-RBC method and its extension to systems with moving loads. The two-level PR-RBC method can be applied to any linear time-domain PDE which admits an affine representation of the parameter. The latter can be recovered by means of empirical quadrature procedure (EQP) \cite{PateraYano} or empirical interpolation method EIM \cite{Barraultetal,Grepletal} if needed. Within our context of SBC, we restrict ourselves to the PDE of linear elastodynamics in this work.

\subsection{Time-Domain Equation}
\label{ss:time_eq}

Let $\Omega\subset\RR^d$, $d=2$ be a bounded domain, with boundary $\partial\Omega$. The proposed method can be naturally applied to the 3D-elastodynamics case, but for purposes of presentation we limit ourselves to the 2D case in this work. The boundary $\partial\Omega$ is assumed to be partitioned into $\Gamma^D$ and $\Gamma^N$, such that Dirichlet boundary conditions are imposed on $\Gamma^D$, while natural boundary conditions are satisfied on $\Gamma^N$. Without loss of generality, the Dirichlet boundary conditions are assumed to be homogeneous, and the non-essential boundary conditions are assumed to be of Neumann type; more details on treating non-homogeneous Dirichlet boundary conditions can be found in \cite{BhouriThesis}. The simulation time interval is noted $[0,T_{\rm final}]$, $T_{\rm final}>0$. Let ${X}\equiv\{v\in [H^1({\Omega})]^{d} \ | \ v|_{{\Gamma}^{D}}=0\}$ be the Hilbert space of admissible real-valued functions; ${X}$ is imbued with inner product $(w,v)_{{X}} \equiv \int_{{\Omega}} \nabla w \cdot \nabla v + w v \, dV$ and induced norm $\| w \|_{{X}} \equiv \sqrt{(w,w)_{{X}}}$. The problem parameterization is denoted ${\mu}\in{\calP}$, where ${\calP}\subset\RR^{{n}_P}$ is a suitable compact set.

Let $m(\cdot,\cdot;\cdot) ; c(\cdot,\cdot;\cdot) ; a(\cdot,\cdot;\cdot) :\big(H^1({\Omega})\big)^d\times \big(H^1({\Omega})\big)^d\times\calP\rightarrow\mathbb{R}$ be the bilinear forms corresponding to the mass, damping and stiffness terms of the 2D-elastodynamics equation respectively. These bilinear forms are defined as follows:
\begin{equation}
m(w,v;\mu)\equiv\rho\;\int_{\Omega}w\cdot  v \ dx \ ,
\end{equation}
\begin{equation}
a(w,v;\mu)\equiv\frac{\nu\; E}{(1+\nu)\;(1-2\nu)}\;\int_{\Omega}\frac{\partial w_i}{\partial x_j}\frac{\partial {v_k}}{\partial x_l}\delta_{ik}\delta_{jl} \ dx + \frac{E}{2(1+\nu)}\;\int_{\Omega}\frac{\partial w_i}{\partial x_j}\frac{\partial {v_k}}{\partial x_l}(\delta_{ik}\delta_{jl}+\delta_{il}\delta_{jk}) \ dx \ , \label{eq:a_bil_form_time_domain}
\end{equation}
\begin{equation}
c(w,v;\mu)\equiv\alpha_{\rm Ray}\; m(w,v;\mu) + \beta_{\rm Ray}\; a(w,v;\mu) \ , \label{eq:c_bil_form_time_domain}
\end{equation}
\noindent where in equation (\ref{eq:a_bil_form_time_domain}) we use the convention of summation over repeated indices, $\nu$ denotes the Poisson ratio, $\rho>0$ is the material density, $\alpha_{\rm Ray}>0$ and $\beta_{\rm Ray}>0$ are the Rayleigh damping coefficients, and $E>0$ is the Young's modulus. Here $f(\cdot,\cdot;\cdot) :\big(H^1({\Omega})\big)^d\times[0,T_{\rm final}]\times\calP\rightarrow\mathbb{R}$ refers to linear form corresponding to the Neumann boundary conditions imposed on $\Gamma^N$. In order to realize an efficient offline-online decomposition, the bilinear and linear forms are assumed to have an affine dependence on the parameter $\mu$. In this work, we are interested in moving loads with time-independent shapes. Hence, we recover the affine dependence of the linear form $f(\cdot,\cdot,\mu_f)$ on the parameters $\mu_f$ by means of EIM \cite{Barraultetal,Grepletal}, where $\mu_f$ denote the parameters within $\mu$ that govern $f(\cdot,\cdot,\mu)$. For the EIM training, we consider the location of the load as a parameter instead of its speed. Such parameter spans the whole domain of possible locations of the moving load. As a consequence, the speed of the moving load, noted $V$, is simply replaced by the location of the load, noted $l$, within the parameters $\mu_f$ to form the parameters $\mufreq_f$ that we consider when performing the EIM.
 %We further assume that the linear form $f(\cdot,\cdot;\cdot)$ satisfies the following space-time separation of variable:
%\begin{equation}
%\label{eq:fxft}
%f(v,t;\mu) = f_t(t;\mu_t)\; f_x(v;{\mu_0}) \ , \forall v\in X \ , \forall t\in[0,T_{\rm final}] \ , \forall\mu\in\calP \ ,
%\end{equation}
%\noindent where $\mu_t$ denotes the parameters governing $f_t(t,\cdot)$, and ${\mu_0}$ refers to all remaining parameters, with $\mu=(\mu_t,\mu_0)$.

The variational formulation of the elastodynamics equation then reads as follows: Find $u(t\in[0,T_{\rm final}];\mu)$ such that $\forall t\in[0,T_{\rm final}]$ ,
\begin{equation}
m\Big(\frac{\partial^2 u(t;\mu)}{\partial t^2},v;\mu\Big)+c\Big(\frac{\partial u(t;\mu)}{\partial t},v;\mu\Big)+a\Big(u(t;\mu),v;\mu\Big)=f(v,t;\mu) \ , \forall v\in X \ , \forall t\in[0,T_{\rm final}]\ , \label{eq:VarFormGlobElast}
\end{equation}
\begin{equation}
u(t;\mu)=0 \ , {\rm on} \ \Gamma^D \ , \forall t\in[0,T_{\rm final}] \ ,
\label{eq:BCDirichletGlob}
\end{equation}
\noindent and
\begin{equation}
u(t=0;\mu)=0 \ ; \frac{\partial u}{\partial t}\Big(t=0;\mu\Big)=0 \ . 
\end{equation}

\noindent Extension to non-zero initial conditions can also be considered.

\subsection{Frequency-Domain Equation}
\label{subsec:FreqEq}

Let $\hat X\equiv\{v \ | \ v = w + \ii y \ , w,y\in X\}$ be the Hilbert space of admissible complex-valued functions; ${\hat X}$ is imbued with inner product $(w,v)_{\hat{X}} \equiv \int_{{\Omega}} \nabla w \cdot \overline{\nabla v} + w \bar v \, dV$ and induced norm $\| w \|_{\hat{X}} \equiv \sqrt{(w,w)_{\hat{X}}}$, where $\overline{\ \cdot \ }$ refers to the complex conjugate operator. The complex problem parameterization is denoted $\mufreq=(\mu,\omega)\in\Pfreq$, which corresponds to the real-valued problem parameter concatenated with the angular frequency $\omega$ as an additional parameter; here $\Pfreq\in\RR^{{n}_P+1}$ refers to the augmented compact parameter set.

We assume that we have non-zero damping such that $\int\limits_0^{\infty}|u(t;\mu)|dt<\infty$, and write $u(t;\mu)=\Re \{\ufreq\big[ \mufreq=(\mu,\omega) \big] e^{ \ii \omega t} \}$, where $\Re$ refers to real part. It follows that $\ufreq(\mufreq)$ satisfies the variational formulation of the Helmholtz equation: Find $\ufreq(\mufreq)\in\hat X$ such that:
\begin{equation}
\afreq(\ufreq(\mufreq),v;\mufreq)=\ffreq(v;\mufreq) \ , \forall v\in\hat X \ , \label{eq:VarFormGlobHom}
\end{equation}
\noindent where $\afreq(\cdot,\cdot;\cdot):\hat X\times \hat X\times\Pfreq\rightarrow\mathbb{C}$ is a sesquilinear form given by
\begin{equation}
\afreq(\cdot,\cdot;\cdot)=-\omega^2\; \hat m(\cdot,\cdot;\cdot)+\ii\omega\; \hat c(\cdot,\cdot;\cdot)+ \hat k(\cdot,\cdot;\cdot) \ .
\end{equation}
\noindent Here $\hat m(\cdot,\cdot;\cdot)$, $\hat c(\cdot,\cdot;\cdot)$ and $\hat k(\cdot,\cdot;\cdot)$ are defined as
\begin{equation}
\hat m(w,v;\mufreq)=\rho\;\int_{\Omega}w\cdot \bar v \ dx \ ,
\end{equation}
\begin{equation}
\hat k(w,v;\mufreq)=\frac{\nu\; E}{(1+\nu)\;(1-2\nu)}\;\int_{\Omega}\frac{\partial w_i}{\partial x_j}\frac{\partial \overline{v_k}}{\partial x_l}\delta_{ik}\delta_{jl} \ dx + \frac{E}{2(1+\nu)}\;\int_{\Omega}\frac{\partial w_i}{\partial x_j}\frac{\partial \overline{v_k}}{\partial x_l}(\delta_{ik}\delta_{jl}+\delta_{il}\delta_{jk}) \ dx \ , \label{eq:k_bil_form_freq_domain}
\end{equation}
\begin{equation}
\hat c(w,v;\mufreq)=\alpha_{\rm Ray}\; \hat m(w,v;\mufreq) + \beta_{\rm Ray}\; \hat a(w,v;\mufreq) \ ,
\end{equation}
\noindent where in equation (\ref{eq:k_bil_form_freq_domain}) we use the convention of summation over repeated indices and $ \ffreq(\cdot;\cdot):\hat X\times\Pfreq\rightarrow\mathbb{R}$ is a continuous anti-linear form associated to the Laplace transform of the $f(\cdot,\cdot;\cdot)$ term. 

Since $m(\cdot,\cdot;\mu)$, $c(\cdot,\cdot;\mu)$, $a(\cdot,\cdot;\mu)$ and $f(\cdot,\cdot;\mu)$ are assumed to have an affine dependence on the parameter $\mu$, it follows that $\afreq(\cdot,\cdot;\mufreq)$ and $\ffreq(\cdot;\mufreq)$ also have an affine dependence on the parameter $\mufreq$. 

\subsection{Finite Element Approximation}
\label{sec:FEapp}

\subsubsection{Finite Element Discretization}

In order to approximate equations (\ref{eq:VarFormGlobElast}) and (\ref{eq:VarFormGlobHom}), we consider a suitably refined finite element (FE) Galerkin approximation: a triangulation ${\calT}^h$ for domain $\Omega$; associated conforming FE approximation spaces $X_h^0\subset \big(H^1(\Omega)\big)^d$ for real-valued functions and $\hat X_h^0$ for complex-valued functions, both of dimension $\calN^0_h$. Let $X_h\equiv X\cap X_h^0$ and $\hat X_h\equiv \hat X\cap \hat X_h^0$, and let $\{\varphi_j\}_{j=1,\ldots,\calN_h}$ denote the associated (real) standard FE nodal basis; here $\calN_{h}$ is the dimension of $X_h$ and $\hat X_h$.

\paragraph{Time-Domain equation}
\label{sssec:time}

The FE approximation $\ufe(t\in[0,T_{\rm final}];\mu)$ to $u(t\in[0,T_{\rm final}];\mu)$ can be obtained by projecting equation (\ref{eq:VarFormGlobElast}) on $X_h$: $\ufe(t;\mu)\in X_{h}\ , \forall t\in[0,T_{\rm final}]$ satisfies 
\begin{equation}
m\Big(\frac{\partial^2 \ufe(t;\mu)}{\partial t^2},v;\mu\Big)+c\Big(\frac{\partial \ufe(t;\mu)}{\partial t},v;\mu\Big)+a\Big(\ufe(t;\mu),v;\mu\Big)=
f(v,t;\mu) \ , \forall v\in X_h \ , \forall t\in[0,T_{\rm final}]\ , \label{eq:TimeVariationalFEGlobal}
\end{equation}
\begin{equation}
\ufe(t=0;\mu)=0 \ ; \frac{\partial \ufe}{\partial t}\Big(t=0;\mu\Big)=0 \ .
\end{equation}

\paragraph{Frequency-Domain equation}
\label{sssec:freq}

Similarly to the time-domain equation, the FE approximation $\ufef(\mufreq)$ to $\ufreq(\mufreq)$ can be obtained by projecting equation (\ref{eq:VarFormGlobHom}) on $X_h$: $\ufef(\mufreq)\in \hat X_{h}$ satisfies
\begin{equation}
\afreq(\ufef(\mufreq),v;\mufreq)=\ffreq(v;\mufreq) \  , \forall v\in \hat X_h \ .\label{eq:FreqVariationalFEGlobal}
\end{equation}

\noindent We now present the discrete equations. Let $\Ahattens\in\RR^{\calN_h\times\calN_h}$ and $\underline{\ffreq_h}\in\RR^{\calN_h}$ be the frequency-domain FE matrix and vector defined as %, $\Xtens_{\rm norm}\in\RR^{\calN^h\times\calN^h}$
\begin{equation}
\Big(\Ahattens(\mufreq)\Big)_{qq'}\equiv \afreq(\varphi_{q'},\varphi_q;\mufreq) \ , 1\leq q,q'\leq\calN_h\ , \label{eq:AMatFreq}
\end{equation}
\begin{equation}
\Big(\underline{\ffreq_h}(\mufreq)\Big)_{q}\equiv \ffreq(\varphi_{q};\mufreq) \ , 1\leq q\leq\calN_h\ \label{eq:FvectFreq}
\end{equation}
\noindent The FE basis vector for $\ufef(\mufreq)\in\RR^{\calN_h}$ is then given by
\begin{equation}
\Ahattens(\mufreq)\;\underline{\ufef}(\mufreq)=\underline{\ffreq_h}(\mufreq) \ .
\end{equation}
\noindent Note that $\Ahattens(\mufreq)$ is typically large but very sparse.

\subsubsection{Finite Element -- Finite Difference Discretization}
\label{subsec:fefddisc}

In order to solve equation (\ref{eq:TimeVariationalFEGlobal}), a finite-difference discretization scheme for time marching with $N_t$ time steps is considered. Let $\Delta t = T_{\rm final}/N_t$ ; $t^j\equiv j\; \Delta t, 0\leq j\leq N_t$. Note that the proposed method can be applied to any finite-difference scheme. A particular scheme is selected in this work for sake of clarity: the unconditionally stable Newmark-$\beta$ scheme with $\beta_t=\frac{1}{4}$ and $\gamma_t=\frac{1}{2}$ (such that the average constant acceleration scheme, or mid-point rule, is obtained) \cite{ChibaKako,Nickelletal}.

Let $\ufefd^j(\mu)$, $0\leq j\leq N_t$ denote the finite element -- finite difference solution at time step $t^j$, and $\dufefd^j(\mu)$ and $\ddufefd^j(\mu)$ the corresponding first and second derivatives in time respectively. Since $\ufefd^0(\mu)=0$ and $\dufefd^0(\mu)=0$, then $\ddufefd^0$ is determined as the solution to
\begin{equation}
m\Big(\ddufefd^0(\mu),v;\mu\Big)=f(v,t=0;\mu) \ , \forall v\in X_h \ .\label{eq:FEFDVariationalddu0}
\end{equation}

\noindent The fields $\ddufefd^j(\mu)$, $\dufefd^j(\mu)$ and $\ufefd^j(\mu)$, $1\leq j\leq N_t$, are then determined as the solutions to the following equations, respectively:
\begin{multline}
m\Big(\ddufefd^j(\mu),v;\mu\Big)+\Delta t\gamma_t \ c\Big(\ddufefd^j(\mu),v;\mu\Big) + \Delta t^2\beta_t \ a\Big(\ddufefd^j(\mu),v;\mu\Big)= f(v,t^j;\mu) \\
-c\Big(\dufefd^{j-1}(\mu)+\Delta t(1-\gamma_t) \ \ddufefd^{j-1}(\mu),v;\mu\Big)-a\Big(\ufefd^{j-1}(\mu)+\Delta t \ \dufefd^{j-1}(\mu)+\Delta t^2(1-\beta_t) \ \ddufefd^{j-1}(\mu),v;\mu\Big) \ , \forall v\in X_h \ ,
\end{multline}
\begin{equation}
\dufefd^j(\mu)=\dufefd^{j-1}(\mu)+\Delta t \ \Big[ (1-\gamma_t) \ \ddufefd^{j-1}(\mu)+\gamma_t \ \ddufefd^j(\mu)\Big]
\end{equation}
\begin{equation}
\ufefd^j(\mu)=\ufefd^{j-1}(\mu)+\Delta t \ \dufefd^{j-1}(\mu)+\Delta t^2 \ \Big[\Big(\frac{1}{2}-\beta_t\Big) \ \ddufefd^{j-1}(\mu)+\beta_t \ \ddufefd^j(\mu) \Big] \ . \label{eq:FEFDVariationaluj}
\end{equation}
 
\noindent We note that the scheme is implicit. 

\noindent We next present the discrete equations in matrix form. Let $\Mtens\in\RR^{\calN_h\times\calN_h}$, $\Ctens\in\RR^{\calN_h\times\calN_h}$ and $\Atens\in\RR^{\calN_h\times\calN_h}$ be the mass, damping and stiffness FE matrices, respectively: %, $\Xtens_{\rm norm}\in\RR^{\calN^h\times\calN^h}$
\begin{equation}
\Big(\Mtens(\mu)\Big)_{qq'}\equiv m(\varphi_{q'},\varphi_q;\mu) \ , 1\leq q,q'\leq\calN_h\ , \label{eq:MMat}
\end{equation}
\begin{equation}
\Big(\Ctens(\mu)\Big)_{qq'}\equiv c(\varphi_{q'},\varphi_q;\mu) \ , 1\leq q,q'\leq\calN_h\ , \label{eq:CMat}
\end{equation}
\begin{equation}
\Big(\Atens(\mu)\Big)_{qq'}\equiv a(\varphi_{q'},\varphi_q;\mu) \ , 1\leq q,q'\leq\calN_h\ . \label{eq:AMat}
\end{equation}

\noindent Furthermore, let $\underline{f_{h}^j}\in\RR^{\calN_h}$, $0\leq j\leq N_t$, be the FE vectors corresponding to the linear form of the time-domain variational formulation at time instance $t^j$, and let $\Ttens\in\RR^{\calN_h\times\calN_h}$ be the time marching matrix:
\begin{equation}
\Big(\underline{f_{h}^j}(\mu)\Big)_{q}\equiv f(\varphi_{q},t_j;\mufreq) \ , 1\leq q\leq\calN_h\ ;\label{eq:Fvect}
\end{equation}
\begin{equation}
\Ttens(\mu)=\Mtens(\mu)+\Delta t\;\gamma_t\;\Ctens(\mu)+\Delta t^2\;\beta_t\;\Atens(\mu)\ .
\end{equation}

\noindent Finally, if $\underline{\ufefd^j}(\mu)\in\RR^{\calN_{h}}$, $\underline{\dufefd^j}(\mu)\in\RR^{\calN_{h}}$ and $\underline{\ddufefd^j}(\mu)\in\RR^{\calN_{h}}$, $0\leq j\leq N_t$, denote the FE basis vectors for ${\ufefd^j(\mu)}$, ${\dufefd^j(\mu)}$ and ${\ddufefd^j(\mu)}$ respectively, we initialize $\underline{\ufefd^0}(\mu)=\underline{0}$, $\underline{\dufefd^0}(\mu)=\underline{0}$, and $\underline{\dufefd^0}(\mu)$ solution of
\begin{equation} 
\Mtens(\mu)\;\underline{\ddufefd^0}(\mu) =\underline{f_{h}^0}(\mu) \ ;
\end{equation}
\noindent we then solve for $\underline{\ufefd^j}(\mu)$, $\underline{\dufefd^j}(\mu)$, and $\underline{\ddufefd^j}(\mu)$, for $1\leq j\leq N_t$, from
\begin{equation}
\Ttens(\mu)\;\underline{\ddufefd^j}(\mu) = \Big[ \underline{f_{h}^j}(\mu) - \Ctens(\mu)\; \Big( \underline{\dufefd^{j-1}}(\mu)+\Delta t(1-\gamma_t) \ \underline{\ddufefd^{j-1}}(\mu) \Big) - \Atens(\mu)\;\Big( \underline{\ufefd^{j-1}}(\mu)+\Delta t \ \underline{\dufefd^{j-1}}(\mu) + \Delta t^2(1-\beta_t) \ \underline{\ddufefd^{j-1}}(\mu) \Big) \Big]\ , 
\end{equation}
\begin{equation}
\underline{\dufefd^{j}}(\mu)=\underline{\dufefd^{j-1}}(\mu)+\Delta t \ \Big[(1-\gamma_t) \ \underline{\ddufefd^{j-1}}(\mu)+\gamma_t \ \underline{\ddufefd^{j}}(\mu) \Big]\ , 
\end{equation}
\begin{equation}
\underline{\ufefd^{j}}(\mu)=\underline{\ufefd^{j-1}}(\mu)+\Delta t \ \underline{\dufefd^{j-1}}(\mu)+\Delta t^2 \ \Big[ \Big(\frac{1}{2}-\beta_t\Big) \underline{\ddufefd^{j-1}}(\mu) + \beta_t \ \underline{\ddufefd^{j}}(\mu) \Big]\ .
\end{equation}

\noindent This completes the FE ``truth" discretization.

\subsection{Two-Level PR-RBC Method: Overview}

In this section, we provide a summary of the two-level PR-RBC method and its extension to systems with moving loads. The two-level PR-RBC method is a domain decomposition technique in which the global system is decomposed into smaller components, which are referred to as instantiated components. This decomposition creates an ensemble of parameterized instantiated components which can be mapped to an ensemble of parameterized archetype components; multiple instantiated components of the global system can correspond to the same archetype component. Moreover, the domain decomposition creates an ensemble of ports, defined as the intersection of the closures of each two adjacent instantiated components, with the latter forming a parameterized bi-component system. These ports can also be mapped to an ensemble of reference ports associated with archetype bi-component systems. For simplicity, the ports are presumed to be mutually disjoint, such that the reference port is associated to two local ports. In the two-level PR-RBC method, the port and component parameterizations correspond to the frequency-domain equation, and the corresponding sesquilinear and anti-linear forms are assumed to have an affine dependence on the parameters. By consequence, any global system built from the library of the archetype components will be governed by a frequency-domain PDE whose sesquilinear and anti-linear forms have an affine dependence on the system parameter. The latter is intrinsically related to the parameters considered for the different instantiated components forming the global system.

\subsubsection{Offline Stage}
\label{subsec:PR_RBC_offline}

The two-level PR-RBC offline stage corresponds to the construction of reduced bases to approximate the solution to the frequency-domain equation (\ref{eq:FreqVariationalFEGlobal}) within the parameterized archetype components and over reference ports. Hence, it is informed by the library of archetype components and reference ports, and is independent of any subsequent (feasible) system assembly. For every reference port, a reduced port space is built in order to approximate the solution on the reference port joining each compatible pair of archetype components \cite{SmetanaPatera,Bhouri2020}. 

The PR-RBC offline stage also includes the construction of reduced bubble spaces for port mode liftings, and a reduced bubble space for each archetype component with non-zero linear form \cite{EftangPatera2}. These reduced bases approximate the solution inside the archetype components domains (zero on the ports). Within our context, the archetype component bubble spaces reduce to the bubble spaces for inhomogeneity associated with non-zero source terms. In this work, all reduced bubble spaces (for port mode liftings and for inhomogeneity) are constructed by Proper Orthogonal Decomposition (POD). Finally, in the context of an offline-online decomposition, all parameter-independent sesquilinear and anti-linear forms needed for the PR-RBC online stage (detailed in Section \ref{subsec:PR_RBC_online}) are computed and stored once in an offline stage.

We emphasize the important role of components. In general, the components distribute the parameter domain: we reduce a large problem with many global parameters to many small problems each with just a few (local) parameters. Components also permit consideration of very large systems: even in the PR-RBC offline stage, we are required to solve FE problems over at most pairs of components | never the full system. Also, components provide geometry and topology parametric variation. And finally, components permit us to more easily justify the PR-RBC Offline investment: we may amortize the offline effort not only over many queries for any particular global system, but over all possible global systems in our family.

As detailed in the development of the two-level PR-RBC method \cite{Bhouri2020}, the frequency-domain equation will be used in order to build the PR-RBC reduced spaces. Since $f(t;\mu)\rightarrow u(t;\mu)$ is a linear time invariant system, we can simply consider the frequency-domain representation of the load to be independent of the frequency $\omega$ when forming the PR-RBC reduced spaces in the offline stage. The exact time-dependency of the load $f(t;\mu)$ will be used within the online stage when solving the time-domain equation after projecting it on the final reduced space. In this work, we are interested in moving loads with time-independent shapes. By consequence, in the offline stage, the anti-linear form $\ffreq(\cdot;\mufreq)$ introduced in (\ref{eq:VarFormGlobHom}) is actually taken as the linear form corresponding to the actual load, but at one given location within its domain of existence. Hence, the parameters governing $\ffreq(\cdot;\cdot)$ are taken as the parameters $\mufreq_f$ introduced in Section \ref{ss:time_eq}. The parameters $\mufreq_f$ correspond to the parameters $\mu_f$, with the load location $l$ replacing the load speed $V$. For instance, if we consider a moving load with a Gaussian shape, then the load amplitude and its spatial width can be taken as the parameters forming $\mu_f$, along with the load speed. In the offline stage, the load amplitude and its spatial width are kept as parameters governing $\ffreq(\cdot;\cdot)$, along with the load location instead of its speed such that $\mufreq_f$ correspond to the load amplitude, its spatial width and its spatial location.

\subsubsection{Online Stage}
\label{subsec:PR_RBC_online}
For the online stage, we consider a global system characterized by a global parameter $\mu$ and defined as an assembly of instantiated archetype components. The online stage has two levels. 

\paragraph{Level 1}

The first-level reduction consists of evaluating the PR-RBC solutions to the frequency-domain PDE for the global system at well-selected frequencies using the reduced bases constructed at the offline stage. To this end, we consider a sufficiently rich angular frequency $\omega$ set $\Xi_\omega$, of size $n_\omega$. We also consdier a set $\Xi_{l}$ of random locations of the load within the domain of existence of the moving load, as we did for the EIM training to approximate $f(\cdot,\cdot;\cdot)$ (see Section \ref{ss:time_eq}) and for the offline stage training (see Section \ref{subsec:PR_RBC_offline}). $\Xi_{l}$ has of the same size $n_\omega$ as $\Xi_\omega$. We then define the online-train (o-t) dataset $\Xi_{\rm o-t}=\{\mufreq\equiv(\mu,l,\omega) ; \omega\in\Xi_\omega,l\in\Xi_l\}$. Hence $\Xi_{\rm o-t}$ can be expressed as: $\Xi_{\rm o-t}=\{\mufreq_j, 1\leq j\leq n_\omega  \}$, where $\mufreq_j=(\mu,l_j,\omega_j)$ for $\omega_j\in\Xi_\omega$ and $l_j\in\Xi_l$. Let $\sigma_t^{\rm ref}$ refer to a characteristic time of $f(\cdot,\cdot;\cdot)$; then $\Xi_\omega$ is chosen as
\begin{equation}\label{eq:xi_omega}
\Xi_\omega=\{0,d\omega,\ldots,\omega_{\rm max}\} \ , \ d\omega=\frac{1}{\underline{c_\omega}\; {\sigma_t^{\rm ref}}} \ , \ \omega_{\rm max}=\frac{\overline{c_\omega}}{{\sigma_t^{\rm ref}}}\ , % 0.1 4
\end{equation}
\noindent so that $n_\omega=\overline{c_\omega} \; \underline{c_\omega}+1$. Note that for the PR-RBC offline stage described in Section \ref{subsec:PR_RBC_offline}, the training is performed over a frequency set that at least contains $\Xi_\omega$. Nonetheless, the PR-RBC approach is still general in the sense that we can consider any assembly of components forming a feasible global system. For the numerical example considered in this work, the PR-RBC offline stage is performed over a frequency set that exactly matches $\Xi_\omega$.

Hence, for each $\mufreq_j\in\Xi_{\rm o-t}$, using the pre-computed and stored sesquilinear and anti-linear form evaluations, the reduced bubble functions for port mode liftings and for inhomogeneity can be computed at very small computational cost compared to a full FE evaluation. Invoking again the pre-computed and stored sesquilinear and anti-linear forms evaluations, the sparse system Schur complement can be formed and solved. We opt for a Petrov-Galerkin projection to construct the reduced system Schur complement, such that our test space is spanned by the lifted port modes, and not by the ``harmonic" functions which form the (statically condensed) trial space. These solutions define the PR-RBC approximations. A first crucial point related to efficiency is the relatively low dimension of the PR-RBC space. A second crucial point related to efficiency is the sparsity of the PR-RBC basis: the support of a given basis function does not exceed two instantiated components for lifted port modes, and is further restricted to just one instantiated component for reduced bubble spaces. The level 1 reduction gives $n_\omega$ PR-RBC approximations, noted $\upr(\mufreq_j)$ for $\mufreq_j\in\Xi_{\rm o-t}$, to $n_\omega$ FE solutions $\ufef(\mufreq_j)$ for $\mufreq_j\in\Xi_{\rm o-t}$. 

\paragraph{Level 2}

The second-level reduction consists of building a ``final" reduced basis from the PR-RBC approximations $\{\upr(\mufreq_j) \ ,\mufreq_j\in\Xi_{\rm o-t}\}$ computed in Level 1 by performing a Strong Greedy procedure to identify a reduced space $X_{\rm RB}$ of size $\calN^{\rm RB}$. These inexpensive PR-RBC approximations are considered as ``truth" solutions. Let $X_{\rm RB\,i}$ denote the reduced space of size $i$, constructed prior to the $i$-th iteration of the Strong Greedy algorithm; let $\upri(\mufreq_j)\in X_{\rm RB\,i}$ denote the RB approximation to $\upr(\mufreq_j)$ obtained using the reduced space $X_{\rm RB\,i}$. Since the PR-RBC approximations are considered as high-fidelity solutions within the greedy algorithm, the next snapshot we add to the reduced space $X_{\rm RB\,i}$ among the PR-RBC solutions $\{\upr(\mufreq_j) \ ,\mufreq_j\in\Xi_{\rm o-t}\}$ is selected based on the norm of the error: $\big|\big|\upr(\mufreq)-\upri(\mufreq)\big|\big|_{H^1(\Omega)}$, for $\mufreq\in\Xi_{\rm o-t}$. An efficient computation of the errors $\big|\big|\upr(\mufreq)-\upri(\mufreq)\big|\big|_{H^1(\Omega)}$, for $\mufreq\in\Xi_{\rm o-t}$, such that only the updated quantities depending on the new snapshot are evaluated at every iteration of the Strong Greedy, is detailed in \cite{BhouriThesis} . The Strong Greedy algorithm is stopped at the $i$-th iteration if the relative error:
\begin{equation}
\label{eq:err_SG}
\frac{\max\limits_{\mufreq\in\Xi_{\rm o-t}}\big|\big|\upr(\mufreq)-\upri(\mufreq)\big|\big|_{H^1(\Omega)}}{\max\limits_{\mufreq\in\Xi_{\rm o-t}}\big|\big|\upr(\mufreq)-\upro(\mufreq)\big|\big|_{H^1(\Omega)}}
\end{equation}
\noindent is below a certain threshold $\epsilon$ or if $i=\min(n_\omega,M)$, where $M$ is a prefixed maximum size allowed for $X_{\rm RB}$ \cite{Bhouri2020}.

The time-domain elastodynamics PDE is then projected and solved within the ``final" reduced space $X_{\rm RB}$ in standard fashion. The exact time signature of the linear form is only used in the time marching performed using the reduced space $X_{\rm RB}$. Note that the PR-RBC offline stage is conducted prior to execution of the two-level procedure. Then, in the online stage, for the given parameter value $\mu$, both levels of reduction are invoked, and hence both levels must be computationally fast. In contrast, the PR-RBC offline stage is run only once independently of the number of parameter values considered in the evaluation of the time-domain solution, hence, we do not give as much importance to the computational cost of the offline stage.
 
Similarly to the FE approximation $\ufe(\cdot,\mu)$, the two-level reduced basis approximation, denoted $\urb(t\in[0,T_{\rm final}];\mu)$, is obtained by projecting equation (\ref{eq:VarFormGlobElast}) on $X_{\rm RB}$: $\urb(t;\mu)\in X_{\rm RB}\ , \forall t\in[0,T_{\rm final}]$ ,  such that
\begin{equation}
m\Big(\frac{\partial^2 \urb(t;\mu)}{\partial t^2},v;\mu\Big)+c\Big(\frac{\partial \urb(t;\mu)}{\partial t},v;\mu\Big)+a\Big(\urb(t;\mu),v;\mu\Big)=
f(v,t;\mu) \ , \forall v\in X_{\rm RB} \ , \forall t\in[0,T_{\rm final}]\ , \label{eq:TimeVariationalRBGlobal}
\end{equation}
\begin{equation}
\urb(t=0;\mu)=0 \ ; \frac{\partial \urb}{\partial t}\Big(t=0;\mu\Big)=0 \ .
\end{equation}

\noindent We now incorporate a finite-difference scheme with the same time-discretization notations introduced in Section \ref{subsec:fefddisc}. Let $\urbfd^j(\mu)$ for $0\leq j\leq N_t$ denote the two-level reduced basis-finite difference solution at time step $t^j$, and $\durbfd^j$ and $\ddurbfd^j$ the corresponding first and second derivatives in time respectively. We initialize $\urbfd^0(\mu)=0$, $\durbfd^0(\mu)=0$, and $\ddurbfd^0$ solution of
\begin{equation}
m\Big(\ddurbfd^0(\mu),v;\mu\Big)= f(v,t=0;\mu) \ , \forall v\in X_{\rm RB} \ ;\label{eq:RBFDVariationalddu0}
\end{equation}
\noindent we then solve for $\ddurbfd^j(\mu)$, $\durbfd^j(\mu)$, and $\urbfd^j(\mu)$, for $1\leq j\leq N_t$, from
\begin{multline}
m\Big(\ddurbfd^j(\mu),v;\mu\Big)+\Delta t\gamma_t \ c\Big(\ddurbfd^j(\mu),v;\mu\Big) + \Delta t^2\beta_t \ a\Big(\ddurbfd^j(\mu),v;\mu\Big)=f(v,t^j;\mu) \\
-c\Big(\durbfd^{j-1}(\mu)+\Delta t(1-\gamma_t) \ \ddurbfd^{j-1}(\mu),v;\mu\Big)-a\Big(\urbfd^{j-1}(\mu)+\Delta t \ \durbfd^{j-1}(\mu)+\Delta t^2(1-\beta_t) \ \ddurbfd^{j-1}(\mu),v;\mu\Big) \ , \forall v\in X_{\rm RB} \ ,
\end{multline}
\begin{equation}
\durbfd^j(\mu)=\durbfd^{j-1}(\mu)+\Delta t \ \Big[ (1-\gamma_t) \ \ddurbfd^{j-1}(\mu)+\gamma_t \ \ddurbfd^j(\mu)\Big] \ ,
\end{equation}
\begin{equation}
\urbfd^j(\mu)=\urbfd^{j-1}(\mu)+\Delta t \ \durbfd^{j-1}(\mu)+\Delta t^2 \ \Big[\Big(\frac{1}{2}-\beta_t\Big) \ \ddurbfd^{j-1}(\mu)+\beta_t \ \ddurbfd^j(\mu) \Big] \ , \label{eq:RBFDVariationaluj}
\end{equation}
\noindent respectively.
 
We now provide the matrix equations. Let $\MRB\in\RR^{\calN^{\rm RB}\times \calN^{\rm RB}}$, $\CRB\in\RR^{\calN^{\rm RB}\times \calN^{\rm RB}}$ and $\ARB\in\RR^{\calN^{\rm RB}\times \calN^{\rm RB}}$ be the mass, damping and stiffness two-level reduced basis matrices, respectively, 
\begin{equation}
\MRB(\mu)=  \XRB(\mu)^H\;\Mtens(\mu)\;\XRB(\mu), \label{eq:MRBMat}
\end{equation}
\begin{equation}
\CRB(\mu)=  \XRB(\mu)^H\;\Ctens(\mu)\;\XRB(\mu), \label{eq:CRBMat}
\end{equation}
\begin{equation}
\ARB(\mu)=  \XRB(\mu)^H\;\Atens(\mu)\;\XRB(\mu), \label{eq:ARBMat}
\end{equation}
\noindent where $\cdot^H$ denotes the Hermitian transpose operator and $\XRB\in\RR^{\calN_{h}\times \calN^{\rm RB}}$ is the FE representation of $X_{\rm RB}$. The column $j$ of $\XRB$ corresponds to the coefficients of the RB basis function $j$ from $X_{\rm RB}$ as represented by the FE nodal basis $\{\varphi_i\}_{i=1,\ldots,\calN_h}$. Similarly, let $\underline{f_{\rm RB}^j}\in\RR^{\calN^{\rm RB}}$, $0\leq j\leq N_t$, be the two-level reduced basis vectors corresponding to the linear form of the time-domain variational formulation at time instance $t^j$, and $\TRB\in\RR^{\calN^{\rm RB}\times \calN^{\rm RB}}$ the time marching matrix,
\begin{equation}
\underline{f_{\rm RB}^j}(\mu)= \XRB(\mu)^\dagger\; \underline{f_{h}^j}(\mu) \label{eq:FRBvect}
\end{equation}
\begin{equation}
\TRB(\mu)=\MRB(\mu)+\Delta t\;\gamma_t\;\CRB(\mu)+\Delta t^2\;\beta_t\;\ARB(\mu)\ ,
\end{equation}
\noindent respectively. 

Let $\underline{\urbfd^j}(\mu)\in\RR^{\calN^{\rm RB}}$, $\underline{\durbfd^j}(\mu)\in\RR^{\calN^{\rm RB}}$ and $\underline{\ddurbfd^j}(\mu)\in\RR^{\calN^{\rm RB}}$, $0\leq j\leq N_t$, denote the two-level reduced basis vectors of ${\urbfd^j}(\mu)$, ${\durbfd^j}(\mu)$ and ${\ddurbfd^j}(\mu)$ respectively. We initialize $\underline{\urbfd^0}(\mu)=\underline{0}$, $\underline{\durbfd^0}(\mu)=\underline{0}$, and $\underline{\ddurbfd^0}(\mu)$ solution of
\begin{equation}
\MRB(\mu)\;\underline{\ddurbfd^0}(\mu) = \underline{f_{\rm RB}^0}(\mu) \ ;
\end{equation}
\noindent we then proceed to time-march as
\begin{multline}
\TRB(\mu)\;\underline{\ddurbfd^j}(\mu) =\\
\Big[ \underline{f_{\rm RB}^j}(\mu) - \CRB(\mu)\; \Big( \underline{\durbfd^{j-1}}(\mu)+\Delta t(1-\gamma_t) \ \underline{\ddurbfd^{j-1}}(\mu) \Big) - \ARB(\mu)\;\Big( \underline{\urbfd^{j-1}}(\mu)+\Delta t \ \underline{\durbfd^{j-1}}(\mu) + \Delta t^2(1-\beta_t) \ \underline{\ddurbfd^{j-1}}(\mu) \Big) \Big]\ , 
\end{multline}
\begin{equation}
\underline{\durbfd^{j}}(\mu)=\underline{\durbfd^{j-1}}(\mu)+\Delta t \ \Big[(1-\gamma_t) \ \underline{\ddurbfd^{j-1}}(\mu)+\gamma_t \ \underline{\ddurbfd^{j}}(\mu) \Big]\ , 
\end{equation}
\begin{equation}
\underline{\urbfd^{j}}(\mu)=\underline{\urbfd^{j-1}}(\mu)+\Delta t \ \underline{\durbfd^{j-1}}(\mu)+\Delta t^2 \ \Big[ \Big(\frac{1}{2}-\beta_t\Big) \underline{\ddurbfd^{j-1}}(\mu) + \beta_t \ \underline{\ddurbfd^{j}}(\mu) \Big]\ , 
\end{equation}
\noindent for $1\leq j\leq N_t$.

Since $X_{\rm RB}$ is a complex-valued reduced basis, the actual two-step reduced basis approximation of the finite element - finite difference solution $\ufefd^j(\mu)$ is given by:
\begin{equation}
\urbfd^j(\mu)=\sum\limits_{k=1}^{\calN_{h}}\Big(\Re\Big[\XRB(\mu)\;\underline{\urbfd^{j}}(\mu)\Big]\Big)_k\varphi_{k}\ , 1\leq j\leq N_t\ .
\end{equation}

\noindent If a quantity of interest $\underline{q}\in\RR^{N_q}$ is considered, then the corresponding FE output matrix $\Qtens\in\RR^{N_q\times\calN_{h}}$ is multiplied by $\XRB(\mu)$ once (offline) to obtain
\begin{equation}
\QRB(\mu)=  \Qtens\;\XRB(\mu), \label{eq:QRBMat}
\end{equation}
\noindent and the two-level reduced basis approximation of $\underline{q}$ at time step $t^j$ is simply given by:
\begin{equation}
\underline{q^j_{\rm RB}}= \Re\Big[ \QRB(\mu)\;\underline{\urbfd^{j}}(\mu)\Big] \ ; 1\leq j\leq N_t. \label{eq:QRB}
\end{equation}

\section{Simulation-Based Classification}\label{sec2:SBC}

%quam. Sed diam turpis, molestie vitae, placerat a, molestie nec, leo\cite{Rothermel1998} Maecenas lacinia. Nam ipsum ligula, eleifend

%Example for bibliography citations cite\cite{Elbaum2002}, cites\cite{Allen2011,Yoo2007}

In this section, a SBC approach is presented. The method takes advantage of features that are built based on the time-domain correlation functions. The two-level PR-RBC method detailed in section \ref{sec:2levelPRRBC} is used to obtain the temporal response of the system needed to compute the two-point correlation functions and by consequence to inexpensively form quasi-exhaustive synthetic training datasets of the features for the classification task. In order to obtain good classification results, a sufficiently large dataset should be constructed, which will be computed with reasonable computation time thanks the two-level PR-RBC method.

First, the global parameter, defined for a given system, can be regarded as follows:
\begin{equation}
\mu = (  \mud, \mun ,\muc) \in \calP \ , \calP = \calPd \times \calPn \times \calPc \; \label{eq:param}
\end{equation}
% where $\calP\subset\RR^{n_P}$ is a sufficiently large compact parameter set 
\noindent and $\mu$ follows the probability density function $\rho_\mu \equiv (\rhod\,\rhon\,\rhoc)$. $\mud$ refers to the parameters that determine the damage status of the structure. This could be the existence of a crack, or the Young modulus for a subdomain to $\Omega$ to model the loss of stiffness as a damage instance. Note that due to its great flexibility in topology and geometry, the two-level PR-RBC method can naturally accommodate for crack existence and is thus particularly well-suited to SHM. $\muc$ refers to the parameters that the operator can control like the location of the sensors used to monitor the structure or any parameter relative to the operation conditions that can be chosen (e.g. load parameters, engine parameters \ldots) and $\mun$ englobes any other parameter that is neither controlled by the operator nor does it define the damage status. This includes for instance material properties that are determined with uncertainty | such as damping coefficients for elastodynamic problems for example | or deviations relative to the parameters used in the model | like the excitation characteristics for instance.

Let $C^d$ denote the set containing the indices of the components that can present a damage. For $c\in C^d$, the locations of $n_{\rm sensors}(c)$ sensors (noted as $x_i^{\rm sensor}(c) \ , 1\leq i\leq n_{\rm sensors}(c) $) are specified and the PR-RBC-based approximation of the outputs for a parameter $\mu$ is defined as: 
\begin{equation}
\underline{q^{c,j}_{\rm RB}}(\mu) \equiv \begin{bmatrix} \urbx(x_1^{\rm sensor}(c) ,j\;t^j,\mu) \\ \vdots \\ \urbx(x_{n_{\rm sensors}(c) }^{\rm sensor}(c) ,t^j,\mu) \\ \urby(x_1^{\rm sensor}(c) ,t^j,\mu) \\ \vdots \\ \urby(x_{n_{\rm sensors}(c) }^{\rm sensor}(c) ,t^j,\mu)\end{bmatrix} \ , 1 \le j \le N_t \ , \Delta t=\frac{T_{\rm final}}{N_t} \ , c\in C^d \ .
\end{equation}

\noindent where $\urb(\cdot,\cdot,\mu)$ refers to the two-level PR-RBC approximation of the displacement for a parameter $\mu$, and in this case $n_{\rm outputs}(c) =2\times n_{\rm sensors}(c)$.

\subsection{Correlation Functions-Based Features}
\label{sec:corrfct}

Appropriate choice of features is absolutely crucial for classification. Features which are sensitive to the anticipated damage but relatively insensitive to nuisance parameters and measurement noise greatly simplify the classification task and ultimately improve the robustness and hence performance of the deployed classifier. Within those considerations, structural damage detection methods using correlation functions of vibration response under stochastic excitation have been developed \cite{Zhang2014,Yang2009,Huo2016}. The normalized two-point correlation function $C_{i,j,k,l}(\tau,\mu)$is defined as:
%\begin{equation}
%\mu \in \calP \rightarrow \Big\{C_{i,j,k,l}(\tau,\mu)\equiv\frac{1}{T_{\rm final}}\int_0^{T_{sim}-\tau}u_{k}(x_{\rm sensor}^i,t,\mu)u_{l}(x_{\rm sensor}^j,t+\tau,\mu)dt\Big\}_{1 \le i,j \le n_{\rm sensors}, 1\le k,l\le d, 0 \le \tau \le T_{\rm final}} \ ,
%\end{equation}
\begin{multline}
C_{i,j,k,l}(\tau,\mu,c)\equiv {\small \frac{1/T_{\rm final}}{\max\limits_{\substack {1 \le p \le n_{\rm sensors}(c) \\ 0 \le t \le T_{\rm final}} }u_{k}(x_{\rm sensor}^p(c),t,\mu) \max\limits_{\substack {1 \le p \le n_{\rm sensors}(c) \\ 0 \le t \le T_{\rm final}} }u_{l}(x_{\rm sensor}^p(c),t,\mu)} } \times \int_0^{T_{sim}-\tau}u_{k}(x_{\rm sensor}^i(c),t,\mu)u_{l}(x_{\rm sensor}^j(c),t+\tau,\mu)dt \\ , \ 1 \le i,j \le n_{\rm sensors}(c), 1\le k,l\le d, 0 \le \tau \le T_{\rm final}
\end{multline}
\noindent
where $T_{\rm final}$ denotes the simulation time and $u_k(x,\cdot,\cdot)$ refers to the displacement in the $k-$th direction at the point $x$ and $x_{\rm sensor}^i(c)$ gives the location of the $i-$th sensor belonging to the component $c$. Based on the choice of its integrand, the two-point correlation function clearly depends on the wave speed since it counts for two displacements with a shift in time. Therefore, it is expected to be sensitive to anticipated damage such as loss of stiffness or crack existence, since these instances of damage considerably affect the wave speed. Using normalization and exploring different strategies for the choice of the time shift; the sensors; and the direction of the displacement considered for the correlation function, the goal is to build features which are not only sensitive to the damage | already guaranteed by the definition of the correlation function | but also relatively insensitive to nuisance variables and measurement noise. Features based on the correlation function amplitude (CCFA) method \cite{Huo2016}, on the inner product vector (IPV) method \cite{Yang2009} and on the AMV method \cite{Zhang2014} were considered. The best results were obtained using the CCFA and the IPV-based features. Hence, in this work, a particular attention is given to the IPV-based features for sake of clarity and conciseness. For $c\in C^d$ and a global parameter $\mu$, we define the IPV-based feature $\calF^{\rm IPV}_c(\mu)$ as:
\begin{equation}
%(\mu \in \calP,c\in C^d) \xrightarrow{\;\,F_c\circ Q_c\;\,} f_c(\mu) \equiv [C_{i,j,k}(\tau=0,\mu,c))_{1 \le i,j \le n_{\rm sensors}(c), 1\le k\le d}]: {\rm IPV} \ ,
(\mu \in \calP,c\in C^d) \rightarrow \calF^{\rm IPV}_c(\mu) \equiv \Big[\big(C_{i,j,k,k}(\tau=0,\mu,c)\big)_{1 \le i,j \le n_{\rm sensors}(c), 1\le k\le d}\Big] \ ,
\end{equation}
but for which no reference sensor needs to be chosen as done in the work of \cite{Yang2009}. Therefore, for the IPV feature, $\nfea(c)=d\times n_{\rm sensors}(c)^2$. Features containing the two-point correlation function only for one direction can also be considered by defining for instance the following ${\rm IPV_x}$ feature, noted $\calF_c^{\rm IPV_x}(\mu)$ for which $k=1$:
\begin{equation}
%(\mu \in \calP,c\in C^d) \xrightarrow{\;\,F_c\circ Q_c\;\,} f_c(\mu) \equiv [C_{i,j,1}(\tau=0,\mu,c))_{1 \le i,j \le n_{\rm sensors}(c)}]: {\rm IPV_x} \ ,
(\mu \in \calP,c\in C^d) \rightarrow \calF_c^{\rm IPV_x}(\mu) \equiv \Big[\big(C_{i,j,1,1}(\tau=0,\mu,c)\big)_{1 \le i,j \le n_{\rm sensors}(c)}] \ ,
\end{equation}
%and ${\rm IPV_y}$ the features for which $k=2$:
%\begin{equation}
%%(\mu \in \calP,c\in C^d) \xrightarrow{\;\,F_c\circ Q_c\;\,} f_c(\mu) \equiv [C_{i,j,2}(\tau=0,\mu,c))_{1 \le i,j \le n_{\rm sensors}(c)}]: {\rm IPV_y} \ ,
%(\mu \in \calP,c\in C^d) \rightarrow f_c(\mu) \equiv [C_{i,j,2}(\tau=0,\mu,c))_{1 \le i,j \le n_{\rm sensors}(c)}]: {\rm IPV_y} \ ,
%\end{equation}
such that for the ${\rm IPV_x}$ feature, $\nfea(c)=n_{\rm sensors}(c)^2$.

\subsection{Classification Task}
\label{sec:class}

In order to perform the classification task, binary damage classes are defined for every component in $c\in C^d$ by considering the mapping $D_c$ defined as:
\begin{equation}
\mud \in \calPd \xrightarrow{\;\,D_c \,\:} l \in \{1,2\} \ , c\in C^d \ , %\ncla
\end{equation}

\noindent where by convention, $l=1$ corresponds to no damage case for component $c$ and $l=2$ to the damaged case. Let $F_c$ denote the mapping of the outputs $\{\underline{q^{c,j}_{\rm RB}}(\mu)\}_{1\leq j\leq N_t}$ to the feature $\calF_c(\mu)$ for a given system parameter $\mu$ as follows:
\begin{equation}
\{\underline{q^{c,j}_{\rm RB}}(\mu)\}_{1\leq j\leq N_t} \xrightarrow{\;\,F_c\;\,} \calF_c(\mu) \in \RR^{\nfea(c)} \ , c\in C^d \ ,
\end{equation}
\noindent where in our case the feature $\calF_c(\mu)$ will either be $\calF_c^{\rm IPV}(\mu)$ or $\calF_c^{\rm IPV_x}(\mu)$.

Finally, for $c\in C^d$, let $Q^{\rm RB}_c$ denote the two-level PR-RBC-based computation process of the outputs $\{\underline{q^{c,j}_{\rm RB}}(\mu)\}_{1\leq j\leq N_t}$ for a given system parameter $\mu$ such that :
\begin{equation}
\mu \xrightarrow{\;\,Q_c^{\rm RB}\;\,} \{\underline{q^{c,j}_{\rm RB}}(\mu)\in \RR^{n_{\rm outputs}(c)}\}_{1\leq j\leq N_t} \ , c\in C^d \ .
\end{equation}

\noindent Then, the considered classifier $C_c$ is applied:
\begin{equation}
\calF_c(\mu) \in \RR^{\nfea(c)} \xrightarrow{\;\,{C_c}\,\;} {l} \in \{1,2\} \ ,
\end{equation}

\noindent and the classifier that minimizes the expected misclassification is selected: 
\begin{equation}
{C_c}^* \equiv \arg \min_{C_c}\; \EE_{\mu}\Big[ L \Big( D_c(\mud), {C_c}(\,F_c(Q^{\rm RB}_c(\mu))\, )\Big) \Big] \ ,
\end{equation}

\noindent where
\begin{equation}
L(k,k') = {\scriptsize \begin{cases} 0 \quad k = k' \\ 1 \quad k \neq k' \end{cases}}\quad k,k' \in \NN^q , \ {\rm for \ some} \ q\in\NN \ .
\end{equation}

As for the classifiers, Artificial Neural Network (ANN) \cite{Bishop1995} and Support Vector Machine (SVM) \cite{Cristianini2000,Cortes1995} are considered. Off-the-shelf Matlab implementations of One-vs-all SVM and ANN are considered. More precisely, fitcsvm Matlab function (using Sequential Minimal Optimization to solve the dual problem) is used for binary Soft-SVM with 10-fold cross validation and the maximum of the score (likelihood that a label comes from a particular class) is consider as the classification criteria. Matlab train function is called for the ANN with $1$ hidden layer containing $10$ nodes and tan-sigmoid activation function. The output layer has a softmax function and the loss function considered is cross-entropy. A scaled cross-conjugate gradient back-propagation method is used with $80\%$ of the dataset for training and $20\%$ for validation. The training stops if validation error does not decrease after 6 iterations. Finally, the maximum number of epochs is fixed at $1000$.

Once a deployed structure is monitored, the classifiers would be tested on features computed from experimental data and which can be expressed as:
\begin{equation}
\fphy_c = \calF_c(q_c) + \eta + \epsilon_b
\end{equation}
\noindent where $\eta$ is the measurement noise and $\epsilon_b$ the bias introduced by relying on the mathematical model to represent the deployed structure dynamics. Hence deviations from the numerical feature in the test dataset should be considered such that for some system parmaeter $\mu \in \calP$, $\fphy_c$ is close to $\calF_c(Q^{\rm RB}_c(\mu))$. Therefore, classifiers are trained on features computed using noiseless data, and  tested on features based on data perturbed with noise $\eta$, such that for every sensor located at $x^{\rm sensor}_i(c)$ , $c\in C^d$, $1\leq i\leq n_{\rm sensors}(c)$, the displacement $\urbk(x^{\rm sensor}_i(c),j\; \Delta t,\mu)$ for $1\leq k\leq d$ , $1\leq j\leq N_t$ is replaced by:
\begin{equation*}
\urbk(x^{\rm sensor}_i(c),t^j,\mu)+\eta \ ,
\end{equation*}
\noindent where $\eta$ is randomly sampled following the Gaussian probability density function 
\begin{equation*}
\calN(0,\sigma\times\max\limits_{1 \leq j \leq N_t} \{\urbk(x^{\rm sensor}(c)_i,t^j,\mu)\} ) \ ,
\end{equation*}
\noindent where $\sigma$ is a predefined noise factor. 

Then, a feature train-test dataset $\Xitt$ of size $\ntt$ is considered, such that the classifier is trained on the noiseless subset $\Xitt_{\phi}$ of size $\phi\times\ntt$ and tested on the set $\Xitt_{1-\phi}$ containing the remaining points of $\Xitt$ which are not picked in $\Xitt_\phi$. (Hence, $\Xitt_{1-\phi}$ is of size $(1-\phi)\times\ntt$). Moreover, the data within the set $\Xitt_{1-\phi}$ is perturbed with noise $\eta$ using the noise factor $\sigma$ as explained above. Finally, $n_{\rm part}$ different random partitions for each $\Xitt$ are considered and the misclassification error averaged over the different random partitions is reported.

The classification approach described above is conduct for every component in $C^d$ and is summarized in the Train-Test Synthetic Supervised Learning algorithm (Algorithm \ref{alg:TT_learn}). In steps \ref{add_noise1} and \ref{add_noise2}, the noise is added to the test dataset as explained above. 
% \mbox{}\hphantom{fo} 

\begin{algorithm}[H]
\caption{T-T-Learning algorithm}\label{alg:TT_learn}
\begin{algorithmic}[1]
%\Function{Strong Greedy}{$\{\upr(\mufreq_j) \ ,\mufreq_j\in\Xi_{\rm o-t}\},n_\omega,M,\epsilon$}

 \State Create t(rain)-t(est) parameter sample: 

$\Xitt \equiv \{\mutt_j\}_{j = 1,\ldots,\ntt}, \; \mutt_j \text{ i.i.d.} \sim \rho_\mu$
\For {$p:=1 {\rm \ to \ } n_{\rm part}$}
\State Consider train-test partition of $\Xitt\equiv\Xitt_{\phi}\bigcup\Xitt_{1-\phi}$
\For {$c {\rm \ in \ } C^d$}
\State Train Classifier $C_c$ on $\Xitt_{\phi}$ without noise
\For {$\mu\in\Xitt_{1-\phi}$}
\For {$i:=1 {\rm \ to \ } n_{\rm sensors}(c)$}
 \For {$k:=1 {\rm \ to \ } d$}
\For {$j:=1 {\rm \ to \ } N_t$}
\State $\eta \sim \calN(0,\sigma\times \max\limits_{1 \leq j \leq N_t} \{\urbk(x^{\rm sensor}_i(c),t^j,\mu)\} )$ \label{add_noise1}
\State $\urbk(x^{\rm sensor}_i(c),t^j,\mu) \leftarrow \urbk(x^{\rm sensor}_i(c),t^j,\mu)+\eta$ \label{add_noise2}
\EndFor
\State \textbf{end for}
\EndFor
\State \textbf{end for}
\EndFor
\State \textbf{end for}
\EndFor
\State \textbf{end for}
\State Test Classifier $C_c$ on $\Xitt_{1-\phi}$ perturbed by noise $\eta$
\State $\text{Err}_{\eta;c}(p)\leftarrow \frac{1}{(1-\phi)\ntt}\displaystyle \sum_{\mu \in \Xitt_{1-\phi}} L \Big(D_c(\mu), C_c( \,\calF_c(\, Q_c^{\rm RB}(\mu)\,))\Big)$
\EndFor
\State \textbf{end for}
\State $\text{Err}_{\eta;c} \leftarrow \frac{1}{n_{\rm part}} \displaystyle \sum_{p=1}^{n_{\rm part}} \text{Err}_{\eta;c}(p)$
\EndFor
\State \textbf{end for}
\State \textbf{return} $\{\text{Err}_{\eta;c}\}_{c\in C^d}$
\end{algorithmic}
\end{algorithm}

Afterwards, the classification task to estimate the binary state of the whole structure is conducted, such that the structure is considered as damaged if there is at least one damaged component. The corresponding binary damage classes are defined by considering the mapping $D_b$ defined as:

\begin{equation}
\mud \in \calPd \xrightarrow{\;\,D_b \,\:} l \in \{1,2\} \ , %\ncla
\end{equation}
such that:
\begin{equation}
D_b(\mud) =  \begin{cases} 2 \ {\rm if \ there \ is \ any} \ c\in C^d \ {\rm s.t.} \ D_c(\mud)=2 \ , \\ 1 \ {\rm if \ for \ any} \ c\in C^d \ , \ D_c(\mud)=1 \ . \end{cases}
\end{equation}

No additional classification training is required for this task since the binary structure state, noted $C_b(\mu)$, can be inferred using the classifiers $C_c$ for $c\in C^d$ as follows:
\begin{equation}
C_b(\mu) =  \begin{cases} 2 \ {\rm if \ there \ is \ any} \ c\in C^d \ {\rm s.t.} \ {C_c}(\,\calF_c(Q^{\rm RB}_c(\mu))\, )=2 \ , \\ 1 \ {\rm if \ for \ any} \ c\in C^d \ , \ {C_c}(\,\calF_c(Q^{\rm RB}_c(\mu))\, )=1 \ . \end{cases}
\end{equation}
\noindent Then, the expected misclassification for the binary structure state defined as:
\begin{equation}
\text{Err}_{\eta;s}\equiv\EE_{\mu}\Big[ L \Big( D_b(\mud), C_b(\mu)\Big) \Big] \ , \label{eq:misclass_bs}
\end{equation}

\noindent can be estimated. Finally, the classification task to precisely determine the state of all the components that can be damaged is considered. Thus, the goal is to discriminate between the $2^{{\rm card}(C^d)}$ possible values that the vector $[D_c(\mud)]_{c\in C^d}$ can take, where ${\rm card}(C^d)$ refers to the size of the set $C^d$. Hence this classification task tends to precisely determine the exact state of the whole structure based on its possible local damages. The classes for such task can be defined by considering the mapping $D$ defined as:
\begin{equation}
\mud \in \calPd \xrightarrow{\;\,D \,\:} [D_c(\mud)]_{c\in C^d} \ . %\ncla
\end{equation}

No additional classification training is required for this task since the states of all the components $c\in C^d$, noted $C(\mu)$, can be inferred using the classifiers $C_c$ for $c\in C^d$ as follows:
\begin{equation}
C(\mu) =  [{C_c}(\,\calF_c(Q^{\rm RB}_c(\mu))\, )]_{c\in C^d} \ .
\end{equation}
\noindent Then, the expected misclassification for the structure state defined as:
\begin{equation}
\text{Err}_{\eta;a}\equiv\EE_{\mu}\Big[ L \Big( D(\mud), C(\mu)\Big) \Big] \ , \label{eq:misclass_ss}
\end{equation}
\noindent can be estimated. 

To obtain a classifier of sufficient quality, features should be sufficiently discriminating between damage instances and sufficiently insensitive to small perturbations induced by nuisance parameters. The train-test sample size $\ntt$ should also be sufficiently large such that damage cases can be distinguished in the many-query context. Finally, the numerical approximation, corresponding to the two-level PR-RBC solution $\urb(\cdot,\cdot)$ in this work, should be sufficiently accurate and close to the high-fidelity approximation, corresponding to the finite element solution $\ufe(\cdot,\cdot)$ in this work, and thus features do not also significantly deviate from those obtained using the high-fidelity approximation. All these criteria will be verified in the next numerical example.

\section{Bridge Example}
\label{sec:num_ex}
 
In this section, the two-level PR-RBC-based SBC approach is applied to a bridge with a moving $2$-axle vehicle problem. In the model, the vehicle-bridge interaction is considered to be the most significant at the connection between the bridge's decks, since this region generally presents gaps and metallic connectors and thus the vehicle's passing in such region induces vibrations that propagate through the bridge. Therefore, the vehicle-bridge interaction is neglected at other regions than the bridge's decks connection.

\subsection{Problem Definition}
\subsubsection{Archetype Components and Bi-component Systems}

\begin{figure}[H]
\centering
	\includegraphics[scale=0.22]{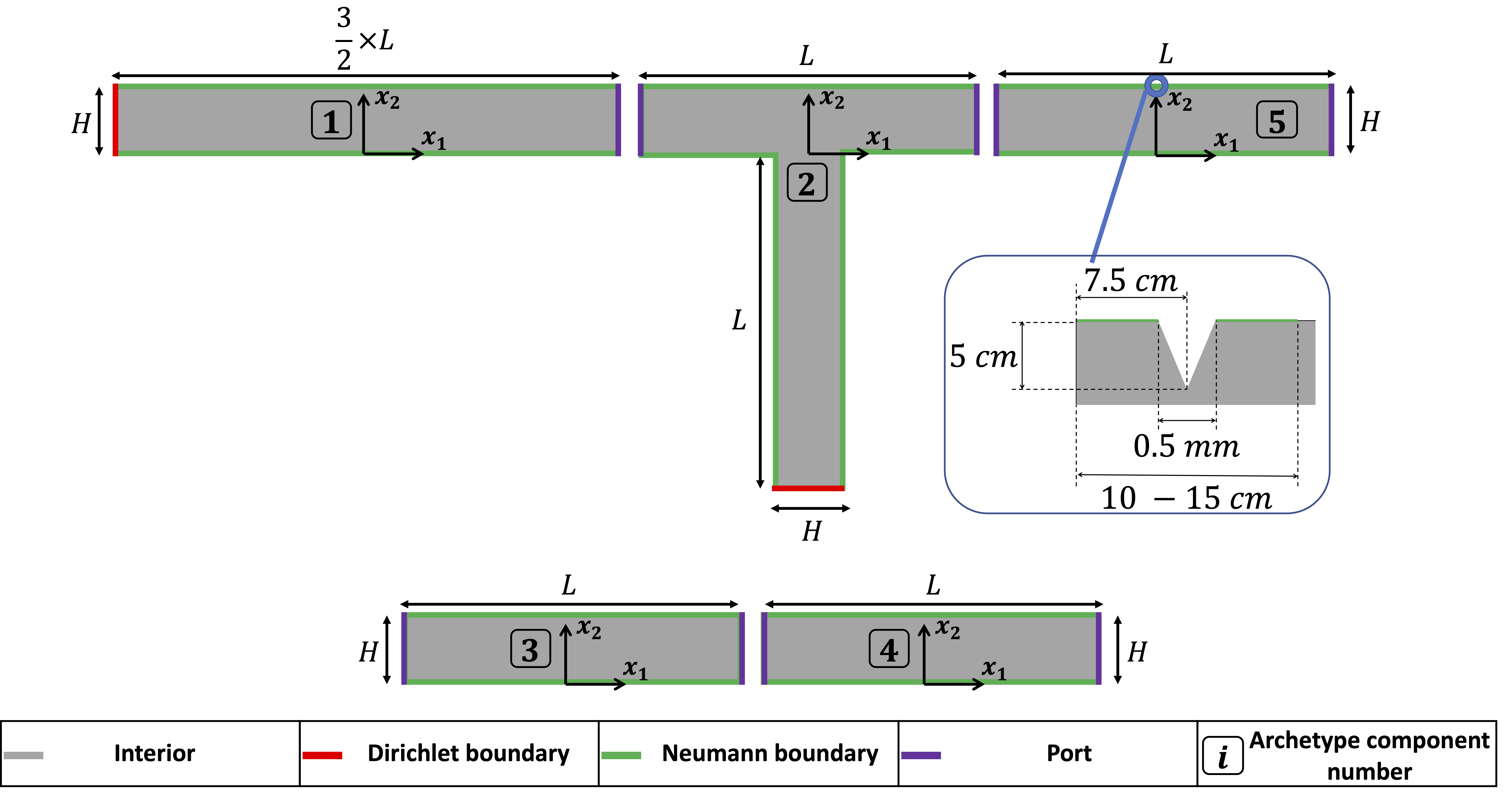}
	\caption{Library of archetype components (considered for construction of elastodynamic bridges with moving vehicle)}					
	\label{fig:bridge_SBC_ref_comp}
\end{figure}

\begin{figure}[H]
	\begin{center}
	\includegraphics[scale=0.22]{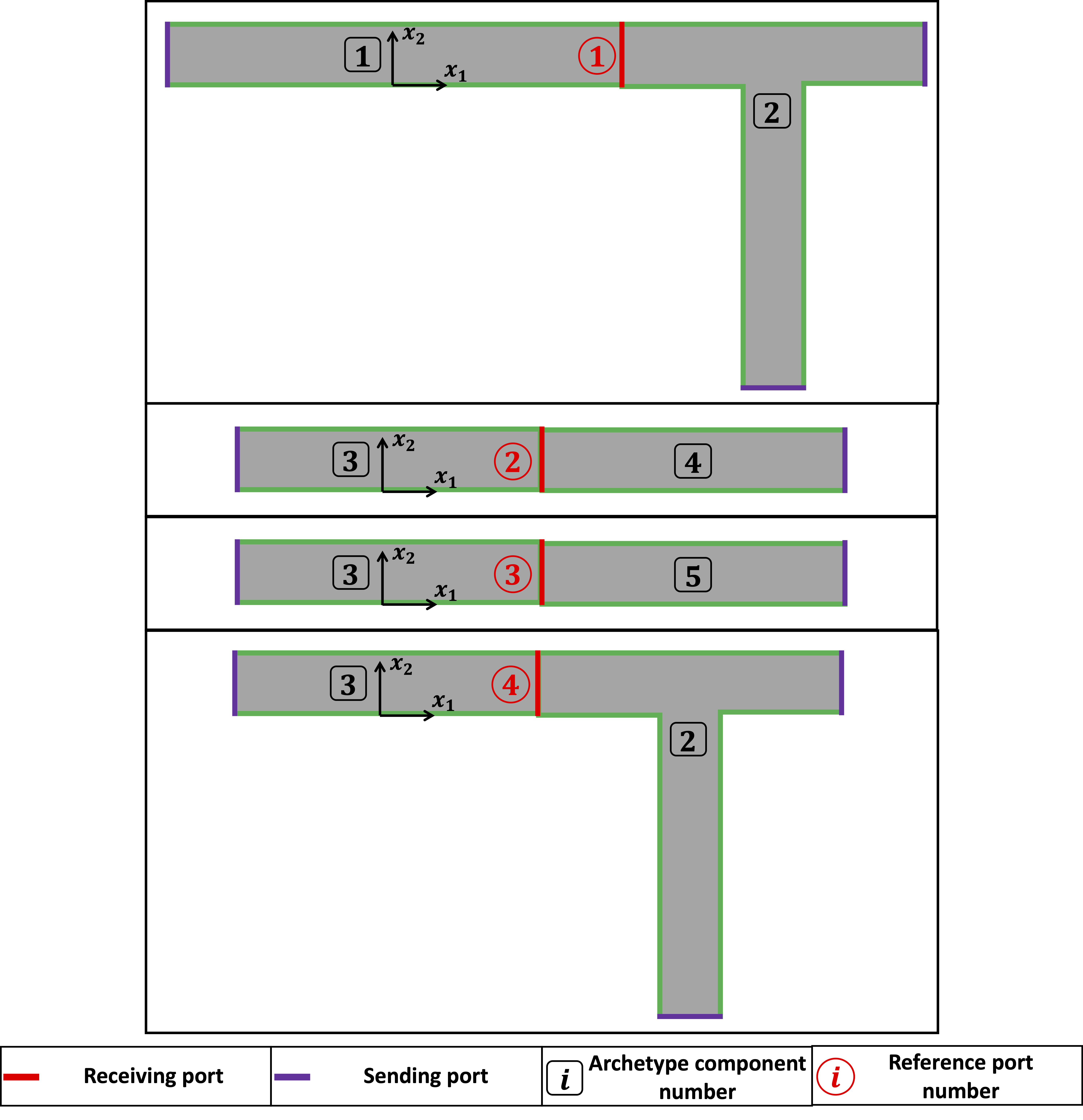}
	\caption{Library of reference ports (archetype bi-components) constructed from the archetype components defined in figure \ref{fig:bridge_SBC_ref_comp}}
	\label{fig:bridge_SBC_ref_port}
	\end{center}	
\end{figure}

To simulate the bridge model described above, the archetype components defined in figure \ref{fig:bridge_SBC_ref_comp} are considered and the associated reference ports are shown in figure \ref{fig:bridge_SBC_ref_port}. The archetype components $1$ has a rectangular geometry of dimension $\frac{3}{2} L\times H$ and contains a homogeneous Dirichlet boundary and homogeneous Neumann boundaries. The archetype component number $2$ has a T shape with thickness equal to $H$ as detailed in figure \ref{fig:bridge_SBC_ref_comp} and also contains a homogenous Dirichlet boundary and homogeneous Neumann boundaries. The archetype components $3$ and $4$ have a rectangular geometry of dimension $L\times H$ and do not contain any Dirichlet boundary. Component $3$ has only homogeneous Neumann boundary conditions, while component $4$ will model the connection between the decks and thus has a non-homogeneous Neumann term corresponding to traction, %boundary condition applied around the midpoint of its upper boundary
\begin{equation}
\sigma(u)\cdot n=\begin{bmatrix} F\; e^{-\frac{[x_1-(-L/2+V\; t)]^2}{\sigma_x^2}}\;\mathbbm{1}_{\{x_2=H\}}\;\mathbbm{1}_{\{-d_1-4\;\sigma_x\leq-L/2+V\; t\leq d_2+4\;\sigma_x\}} \\ -c_{\rm friction}\; F\; e^{-\frac{[x_1-(-L/2+V\; t)]^2}{\sigma_x^2}}\;\mathbbm{1}_{\{x_2=H\}}\;\mathbbm{1}_{\{-d_1-4\;\sigma_x\leq-L/2+V\; t\leq d_2+4\;\sigma_x\}} \end{bmatrix}  \ , \forall x_1\in \Big[-\frac{L}{2},\frac{L}{2}\Big] \ , x_2\in\{0,H\} \ , \label{eq:nh_NBC}
\end{equation}
\noindent applied around the midpoint of its upper boundary. Here $\sigma(u)$ denotes the stress tensor, $F$ the load amplitude, $V$ the vehicle speed, $d_1$ and $d_2$ the limits of the vehicle-bridge interaction, $\sigma_x$ the load spatial width, $c_{\rm friction}$ the fraction coefficient, $\mathbbm{1}_C$ the 2D-function equal to $1$ if the condition $C$ is satisfied and $0$ otherwise, and $n$ the outer normal of the geometric domain. This particular spatial dependence of vehicle-bridge interaction is chosen based on the results obtained in the works of \cite{Yu2017,Yap1989}, which show profiles close to moving Gaussians.

Finally, the archetype component number $5$ has a rectangular geometry of dimension $L\times H$ with non-homogeneous Neumann boundary condition similar to the one applied for the archetype component number $4$ (\ref{eq:nh_NBC}), but archetype component $6$ also contains a crack whose dimensions are detailed in figure \ref{fig:bridge_SBC_ref_comp}. This archetype component $6$ will be used along with the archetype component $4$ to model the existence or not of damage. %From now on, $L_{\rm ref}=H$ and $L$ and $H$ will refer to the dimensionless lengths (which means that $H=1$ but $H$ will not be omitted from the equations for sake of clarity).

Each of archetype components number $1$, $2$ and $3$ has $3$ parameters which consist of the Young's modulus and the two Rayleigh damping coefficients as introduced in equations (\ref{eq:a_bil_form_time_domain}) and (\ref{eq:c_bil_form_time_domain}). Hence, the corresponding frequency-domain variational problem has $4$ parameters | including the angular frequency. Each of archetype components number $4$ and $5$ has $6$ additional parameters defining the non-homogeneous Neumann boundary condition detailed in equation (\ref{eq:nh_NBC}): the load amplitude $F$, the load center defined by the vehicle speed $V$, the load spatial width $\sigma_x$, the friction coefficient $c_{\rm friction}$, and the two parameters $d_1$ and $d_2$ defining the interval length around the midpoint of the upper boundary such that the vehicle-bridge interaction is nonzero if the load center is within the spatial interval $[-d_1-4\;\sigma_x, d_2+4\;\sigma_x]$. The frequency-domain variational problem for archetype components $4$ and $5$ has $10$ parameters in total. The reference ports parameters and the boundary conditions considered for the associated bi-component problems follow naturally from the parameters and the boundary conditions defined for the archetype components. The variational problem for archetype component-pairwise training has either $7$ parameters, for reference ports not involving the archetype component $4$ or $5$, or $13$ parameters for the reference ports involving the archetype component $4$ or $5$. These local parameters associated with the frequency-domain variational problems define the parameters spaces considered in building the PR-RBC reduced bases that are used at Level 1 reduction. 

$L$ and $H$ are taken equal to $5 \ {\rm m} $ and $1  \ {\rm m}$ respectively. The bridge's material is chosen as reinforced concrete since it is one of the most used materials in building bridges. Thus, the Rayleigh damping coefficients $\alpha_{\rm Ray}$ and $\beta_{\rm Ray}$ follow the uniform distributions over the intervals 
\begin{equation*}
[0.566 \ {\rm s^{-1}},4.311 \ {\rm s^{-1}}] \ {\rm and} \ [0.009 \ {\rm s},0.021 \ {\rm s}]
\end{equation*}
\noindent respectively, as established in the work of \cite{Feng2009} for highway bridges in concrete. The Young's modulus also follows a uniform distributions over 
\begin{equation*}
[29 \ {\rm GPa},37 \ {\rm GPa}] % $\bar{E}$ is chosen equal to $\frac{E_{\rm min}+E_{\rm max}}{2}$. 
\end{equation*}
\noindent as determined in the work of \cite{Musial2016} for reinforced concrete. The spatial width of the vehicle-bridge interaction $\sigma_x$ and its amplitude $F$ also follow uniform distributions over the intervals 
\begin{equation*}
[2 \ {\rm cm},4 \ {\rm cm}] \ {\rm and} \ [10^6 \ {\rm Pa},2\times 10^6 \ {\rm Pa}]
\end{equation*}
\noindent respectively, based on the results obtained in the works of \cite{Yu2017,Yap1989}. Given the work of \cite{Yu2017} (and also using https://www.engineeringtoolbox.com), the friction coefficient between a typical car tire and asphalt can be considered following a uniform distribution over \begin{equation*}
[0.5,0.7] \ .
\end{equation*}
\noindent The distance along which the load is applied, and thus the span of the region defining the connection between the bridge's decks, is determined by the two parameters $d_1$ and $d_2$. Based on the work of \cite{Culmo2009} that gives the details of the possible connection options for prefabricated bridge elements in reinforced concrete, $d_1$ and $d_2$ follow a uniform distribution over 
\begin{equation*}
[d_{\rm min}=10 \ {\rm cm}, d_{\rm max}=15 \ {\rm cm}] \ .
\end{equation*}
\noindent Moreover, the geometry of the crack for component $6$ was fixed such that it has the maximum allowable dimensions fixed by Federal Highway Administration of the U.S. Department of Transportation \cite{ACI2001,Balakumaran2018}. Finally. for reinforced concrete, $\rho=2400 \ {\rm kg.m^{-3}}$ and $\nu=0.15$ \cite{Logan2009} (using https://www.engineeringtoolbox.com and https://www.concrete.org/). %We set the Young's modulus mean value: $\bar{E}=\frac{29+37}{2}=33\; GPa$ as nominal value and choose $T_{\rm ref}\equiv\frac{H}{c_t}$, where $c_t=\sqrt{\frac{\bar{E}}{2\;\rho\;(1+\nu)}}$ is the celerity of the transverse wave in infinite domain without damping.

Since the system's excitation is a moving vehicle, to determine the characteristic time of the load $\sigma_t^{\rm ref}$ (see equation (\ref{eq:xi_omega})) and thus the angular frequency set $\Xi_\omega$, the vehicle speed needs to be specified. Vehicles with speed following a uniform distribution between $V_{\rm min}=15 \ {\rm km.h^{-1}}$ and $V_{\rm max}=50 \ {\rm km.h^{-1}}$ are considered. We choose $\sigma_t^{\rm ref}=\frac{d_{\rm min}}{V_{\rm max}}$ and take
\begin{equation}
d\omega=\frac{V_{\rm min}}{2\; d_{\rm max}} \ , \ \omega_{\rm max}=10\;\frac{V_{\rm max}}{2\; d_{\rm min}}\ ,
\end{equation}
\noindent such that $\underline{c_\omega}=10$, $\overline{c_\omega}=5$ and $n_\omega=51$ (see equation (\ref{eq:xi_omega})).

The sizes of the different reduced bases formed at Level 1 reduction are chosen based on the decrease of the eigenvalues of the transfer eigenvalue problem and the decrease of the POD modes for the reduced bubble space for inhomogeneity, the reduced port space, and the reduced space for port mode lifting. Table \ref{tab:bridge_SBC_PR_RBC_RB_sizes} gathers the sizes of the different reduced bases and the computation time to run the offline stage needed by Level 1 reduction. All simulations considered in this work were run on a 4-core laptop (with a 3.5 GHz Intel CPU and 16 GB RAM). We also provide the PR-RBC offline cost, though this cost is amortized over the many online Level 1 - Level 2 queries.

\begin{table}[htb]
\begin{center}\begin{tabular}{|p{10cm}|p{5cm}|} \hline
	Size of training set $\Xi_{\rm o-t}$ & $51$ \\ \hline
	Size of port spaces & \begin{tabular}{@{}l@{}}
                   			$3$ for reference ports  $1$ and $4$ \\
                   			$4$ for reference ports  $2$ and $3$\\
                				\end{tabular} \\ \hline %0.91s
        Size of bubble spaces for port mode lifting & $2$ \\ \hline %0.24s
       Size of bubble space for inhomogeneity for archetype components $4$ and $5$ & $1$  \\ \hline
       Computation time to run PR-RBC offline stage & $64\; s$ \\ \hline %1.37	
\end{tabular}\end{center}
\caption{PR-RBC reduced bases sizes for elastodynamics bridge}
\label{tab:bridge_SBC_PR_RBC_RB_sizes}
\end{table}%

\subsubsection{Global System}

\begin{figure}[H]
\centering
	\includegraphics[scale=0.11]{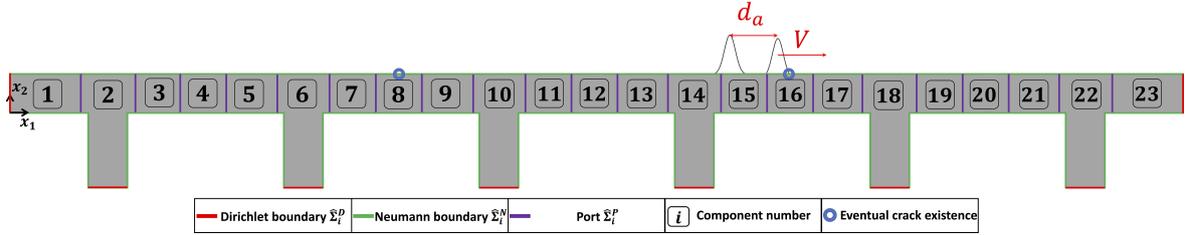}
	\caption{Note-to-scale representation of the global system for the bridge: the Gaussian curves on the top boundary indicate the vehicle's load applied on the bridge}\label{fig:bridge_SBC_glob_model_time}	
\end{figure}

We consider the global system consisting of $n_{\rm comp}=23$ instantiated components presented in figure \ref{fig:bridge_SBC_glob_model_time}. The mapping of each instantiated component of the global system to the corresponding archetype component is given in table \ref{tab:bridge_mapping_comp}. Since each of components number $8$ and $16$ can have a crack, each of them is either mapped to the archetype component $4$ or $5$.

\begin{table}[htb]
\begin{center}\begin{tabular}{|p{7cm}|p{7cm}|} \hline
	Instantiated component number in global system & Archetype component number \\ \hline
	$1,23$ & $1$ \\ \hline 
        $3,5,7,9,11,13,15,17,19,21$ & $3$ \\ \hline
        $2,6,10,14,18,22$ & $2$ \\ \hline
        $4,12,20$ & $4$ \\ \hline
        $8,16$ & $4$ or $5$ \\ \hline  	
\end{tabular}\end{center}
\caption{Instantiated components to archetypes mapping}
\label{tab:bridge_mapping_comp}
\end{table}%

Note that in addition to the parameters defined for every component, two more parameters for the global system have to be considered: the vehicle speed $V$ and the distance between the two axles noted $d_a$ as illustrated in figure \ref{fig:bridge_SBC_glob_model_time}. $d_a$ follows a Gaussian distirubtion with mean equal to $\overline{d_a}=3 \ {\rm m}$ and a standard deviation equal to $\sigma_{d_a}=0.5 \ {\rm m}$ based on generic dimensions of 2-axle vehicles, while $V$ is sampled from a uniform probability density function define over $[V_{\rm min},V_{\rm max}]$. The time-domain variational problem for the global system has a total of $45$ parameters: the bridge's two Rayleigh damping coefficients; the Young's modulus for each component; the vehicle-bridge's interaction width, the load amplitude and the friction coefficient for each axle ($6$ parameters); the vehicle speed $V$; the distance between the two axles $d_a$; the parameters $d_1$ and $d_2$ defining the spatial region along which the load is applied ($10$ parameters); and the two parameters defining the possible existence of crack within components $8$ and $16$. As mentioned above, the frequency-domain variational problem has either $4$ or $10$ parameters for the archetype components, and $7$ or $13$ parameters for the archetype bi-components defining the reference ports. Those sizes are lower than the size of the parameter space of the global domain (equal to $45$) and thus it shows how the two-level PR-RBC method reduces the effective dimensionality of the parameter spaces considered in the variational problems. This reduction of the size of the parameter spaces is even more enhanced as larger global domains with more instantiated components are considered. 

For the time-domain problem approximation, the simulation time $T_{\rm final}$ is fixed based on the vehicle speed $V$ such that it travels across the entire bridge during $T_{\rm final}$. The size of the reduced space $X_{\rm RB}$ constructed by Strong Greedy approach (see Paragraph Level 2 of Section \ref{subsec:PR_RBC_online}) is fixed by imposing $\epsilon=10^{-5}$ as a threshold for the relative error (\ref{eq:err_SG}). As expected for the second order Newmark-$\beta$ scheme ($\beta_t=\frac{1}{4}$, $\gamma_t=\frac{1}{2}$), the order of convergence in time-discretization is equal to $p=2$. The number of time steps $N_t$, and equivalently the step-size $\Delta t$, are fixed based on the second order convergence of the normalized Richardson's extrapolation-based error indicator given by:
\begin{equation}
\epsilon_{\Delta t}\equiv\frac{1}{\max\limits_{1\leq j\leq N_t} \Big|\Big|\urbfd^j(\mu)\Big|\Big|_{H^1(\Omega)}}\;\frac{\max\limits_{1\leq j\leq N_t/2}\Big|\Big|\urbfd^{2\;j}(\mu)-\urbfderr^j(\mu)\Big|\Big|_{H^1(\Omega)}} {2^p-1} \ .
\end{equation}
\noindent The number of time steps $N_t$ is fixed such that we impose $\epsilon_{\Delta t}\leq10^{-4}$. 

For purposes of presentation, we consider four global parameters $\mu_{{\rm example},i}$, $i=1,\ldots,4$ corresponding to the following cases:
\begin{enumerate} \item Case 1: $\mu_{{\rm example},1}$ is such that all parameters are equal to their average value (average of their probability density functions) and we don't have any crack (components $8$ and $16$ are both mapped to archetype component $4$).
\item Case 2: $\mu_{{\rm example},2}$ is such that all parameters are equal to their average value and component $8$ has a crack (component $8$ is thus mapped to the archetype component $5$, while component $16$ is mapped to the archetype component $4$).
\item Case 3: $\mu_{{\rm example},3}$ is such that all parameters, except $d_a$, are equal to their minimum allowable values which consist in the lower bound of the intervals used for the uniform probability density functions and $d_a=\overline{d_a}-4\times\sigma_{d_a}$. Component $16$ has a crack (component $8$ is thus mapped to the archetype component $4$, while component $16$ is mapped to the archetype component $5$).
\item Case 4: $\mu_{{\rm example},4}$ is such that all parameters, except $d_a$, are equal to their maximum allowable values which consist in the upper bound of the intervals used for the uniform probability density functions and $d_a=\overline{d_a}+4\times\sigma_{d_a}$. We also consider 2 cracks (component $8$ and $16$ are both mapped to the archetype component $5$).
\end{enumerate} 
We also consider $10$ randomly sampled global parameters $\mu_{{\rm rand} \, i}$, $i=1,\ldots,10$ and note 
\begin{equation}
\Xi_{\rm o}\equiv\{\mu_{{\rm example},i}, i=1,\ldots,10\} \cup \{\mu_{{\rm rand} \, j}, j=1,\ldots, 10\} \ .
\end{equation}

\begin{figure}[H]
	\begin{center}
	\includegraphics[width=\textwidth]{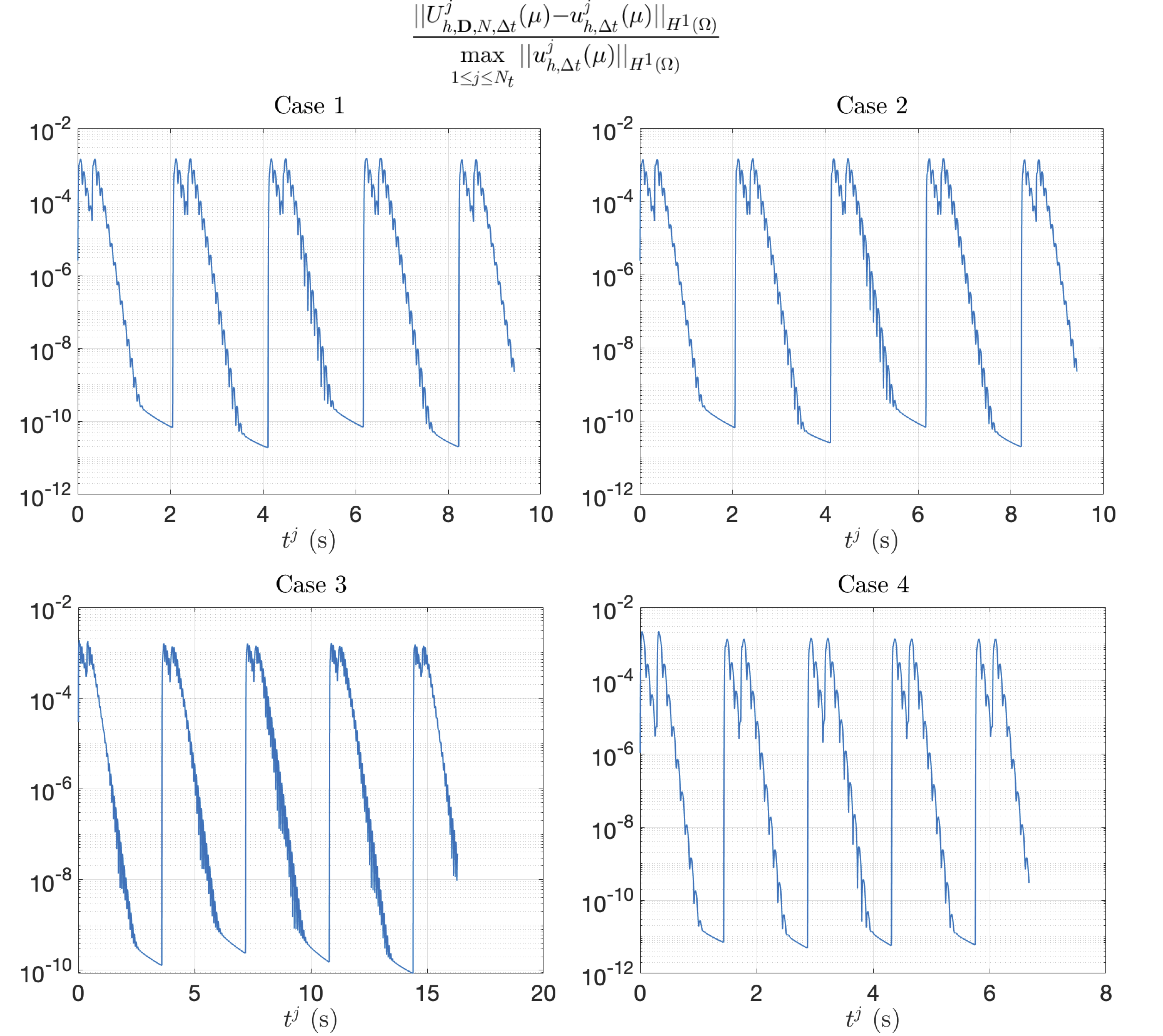}
	\caption{Time-domain relative error for the elastodynamics bridge and the global parameters $\mu_{{\rm example},i}$, $i=1,\ldots,4$}
	\label{fig:err_FE_PR_RBC_by_max_time_bridge_PR_RBC}	
	\end{center}	
\end{figure}

The size of $X_{\rm RB}$ obtained by imposing $\epsilon=10^{-5}$ is on average equal to $N=30$ for the parameters $\Xi_{\rm o}$, while imposing $\epsilon_{\Delta t}\leq10^{-4}$ requires $N_t=10^4$ based on the convergence of $\epsilon_{\Delta t}$ with $N_t$. Table \ref{tab:bridge_SBC_RB_online} gathers the computation time to estimate $\urbfd^j(\mu)$ , $1\leq j\leq N_t$, averaged over the parameters $\Xi_{\rm o}$. The size of the full $\PP_2$ FE approximation space is $\calN^h=17532\pm36$ and computing one full FE simulation with $N_t=10^4$ takes $215\; s$ on average (without performing the Richardson extrapolation). For confirmation purposes, we consider the time-domain relative error between the two-level PR-RBC solution and the FE approximation for a global parameter $\mu$,
\begin{equation*}
\frac{||\urbfd^j(\mu)-\ufefd^j(\mu)||_{H^1(\Omega)}}{\max\limits_{1\leq j\leq N_t} ||\ufefd^j(\mu)||_{H^1(\Omega)}} \ , 1\leq j\leq N_t\ .
\end{equation*}
\noindent Figure \ref{fig:err_FE_PR_RBC_by_max_time_bridge_PR_RBC} gives the evolution of the relative error for the global parameters $\mu_{{\rm example},i}$, $i=1,\ldots,4$ and we can verify that it is well below $1\%$, confirming the sufficiently refined time discretization thanks to the criterion $\epsilon_{\Delta t}<10^{-4}$, and the sufficiently rich reduced space $X_{\rm RB}$ thanks to the strong greedy criterion imposed with $\epsilon=10^{-5}$. 

\begin{table}[htb]
\begin{center}\begin{tabular}{|p{12cm}|p{1cm}|} \hline
	PR-RBC online stage called $n_\omega=51$ times & $2.92 \;s$ \\ \hline
	Strong greedy & $0.56\; s$  \\ \hline %0.91s
        Time marching & $0.25\; s$   \\ \hline % 0.43\; s 0.24s
        Total computation time to estimate $\urbfd^j(\mu)$ , $1\leq j\leq N_t$ & $3.73\; s$ \\ \hline % 3.97\; s 1.37
\end{tabular}\end{center}
\caption{Computation time of the two-level reduced basis method for elastodynamics bridge}
\label{tab:bridge_SBC_RB_online}
\end{table}%

In conclusion, the two-level reduction approach has a computation cost $58$ times lower than the FE simulation for this bridge example. In the next section, we show the usefulness of such computational speedup in the context of Simulation Based Classification for Structural Health Monitoring. 
% $268$ times lower   $11.04$ hours 

\subsection{Classification Task}
\begin{figure}[H]
\centering
	\includegraphics[scale=0.1]{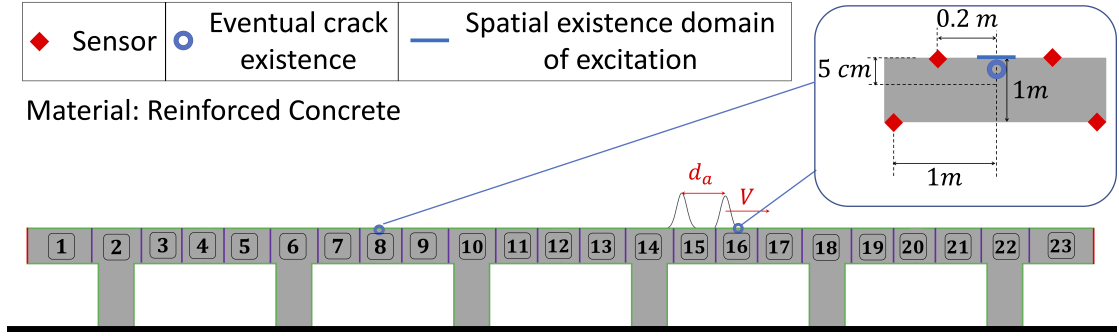}
	\caption{Note-to-scale representation of the global system for the bridge with sensors location}\label{fig:sensors_bridge_SBC}	
\end{figure}

A synthetic train-test dataset is constructed using the two-level PR-RBC method for the SBC task of the global system given in figure \ref{fig:bridge_SBC_glob_model_time}. To do so, the parameters are sampled following their governing distributions introduced in the previous subsection, and the corresponding two-level PR-RBC time-domain solutions are computed. Moreover, for this classification task, components $8$ and $16$ can be damaged based on the existence or not of a crack, such that $C^d=\{8,16\}$. For every component $c\in C^d$, $n_{\rm sensors}(c)=4$ sensors are chosen and located as shown in figure \ref{fig:sensors_bridge_SBC}.

Note that the potential cracks are located at $(8\times L,H)$ and $(16\times L,H)$. The damage is defined as the existence of crack, thus $\mud\equiv(\theta_8,\theta_{16})$, where $\theta_i$ follows a uniform distribution over the set $\{1,2\}$ and component $i$ has a crack if $\theta_i=1$.

The T-T-Learning algorithm (Algorithm \ref{alg:TT_learn}) is run for each component $c\in C^d$ and for different features and classifiers. A particular focus will be given to the effect of the size $\ntt$ of the train-test $\Xitt$ on the accuracy of the classification. For the different results given this section, $\phi$ is taken equal to $0.7$, which means that for every train-test set $\Xitt$ considered, $70\%$ of $\Xitt$ are used to train the classifier, while the remaining $30\%$ are used to test it, and $n_{\rm part}=100$ random partitions are considered.

\begin{figure}[H]
	\centering
	\includegraphics[width=\textwidth]{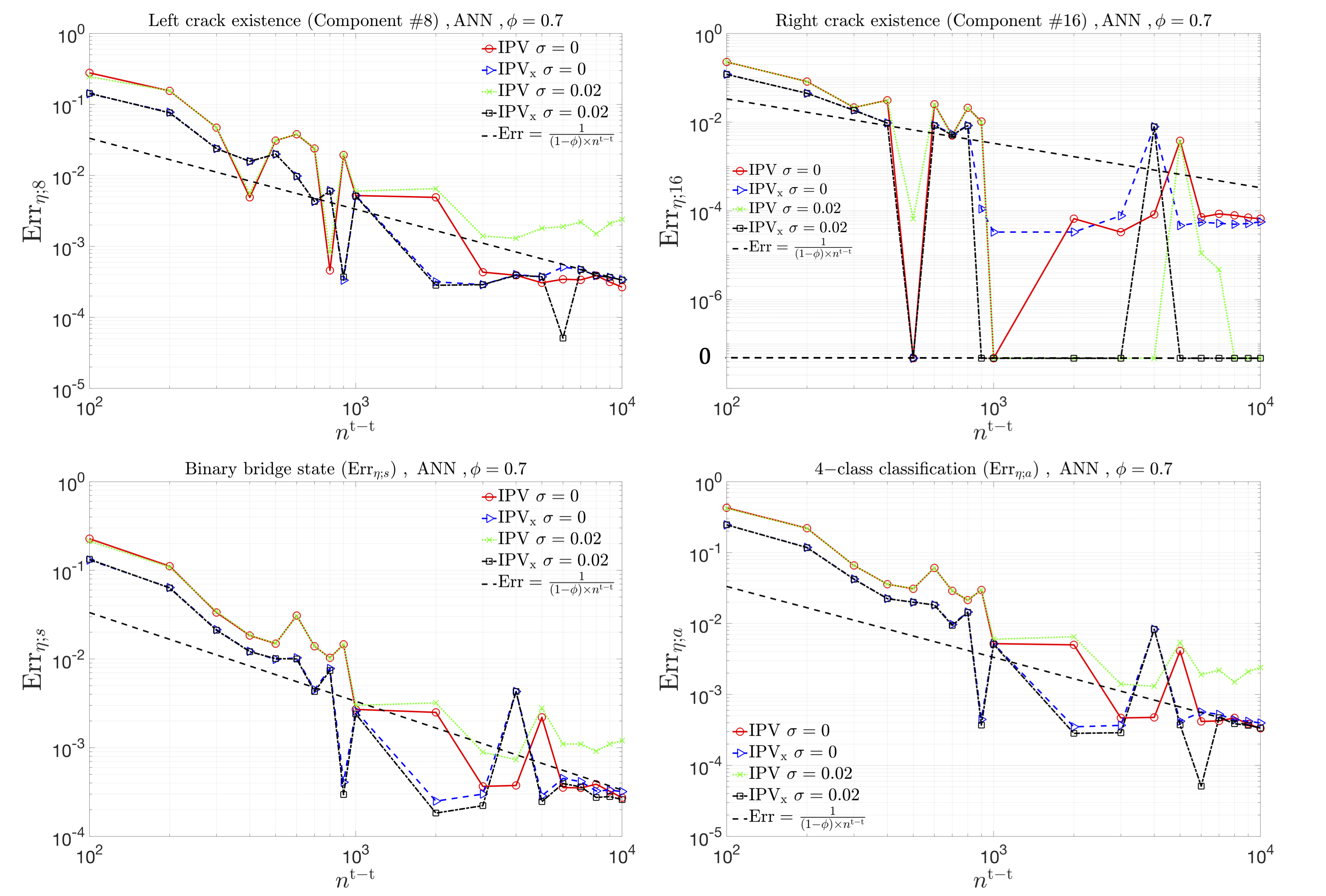} % 0.33
	\caption{Misclassification errors variation with $\ntt$ using ${\rm IPV}$, ${\rm IPV_x}$ features and ANN with $\phi=0.7$}\label{fig:bridge_classification_ANN_new_comp_IPV_IPVx}	
\end{figure}

Figure \ref{fig:bridge_classification_ANN_new_comp_IPV_IPVx} shows the decrease of the expected misclassification for the different components $c\in C^d$, the expected misclassification for the binary structure state and the expected misclassification for the structure state ($4$-class classification task in this example), with the increase of the train-test dataset size $\ntt$ and obtained using the ${\rm IPV}$ and ${\rm IPV_x}$ features and ANN. The figure contains the results obtained for noiseless test dataset and for test dataset perturbed with noise with a factor equal to $\sigma=0.02$. This value was chosen based on error estimates carried out for bridge displacement computed from measured acceleration records \cite{Gindy2007}. Computer vision-based techniques for bridge displacement measurements generally present higher noise levels \cite{Jo2018}. As $\ntt$ gets close to $10^4$, the ${\rm IPV_x}$ feature gives misclassification errors of orders around or below $10^{-4}$. Such observation is still valid even for the structure state classification which is a $4$-class classification task and even with a test dataset perturbed by noise. Figure \ref{fig:bridge_classification_ANN_new_comp_IPV_IPVx} also contains the misclassification error corresponding to only one misclassified point out of the whole test dataset, which is plotted in black dashed line. The misclassification errors obtained correspond to less than one misclassified point as $\ntt$ gets close to $10^4$ for all cases considered except for the ${\rm IPV}$ feature with $\sigma=0.02$. 

As mentioned above, the computation time for one simulation using the two-level PR-RBC approach is equal to $3.73s$, while computing the features takes $0.05s$ (for the $2$ different noise levels considered). All simulations considered in this work were run on a 4-core laptop (with a 3.5 GHz Intel CPU and 16 GB RAM). To obtain the results presented above, a train-test dataset $\Xitt$ of size $\ntt=10^4$ needs to be constructed. Conducting such a task using the two-level PR-RBC approach has a total computation time of $10.52$ hours (taking into account the computational cost of the offline stage, which is run only once), as opposed to an estimated $24.88$ days using the full FE approximation.

\begin{figure}[H]
	\centering
	\includegraphics[scale=0.1]{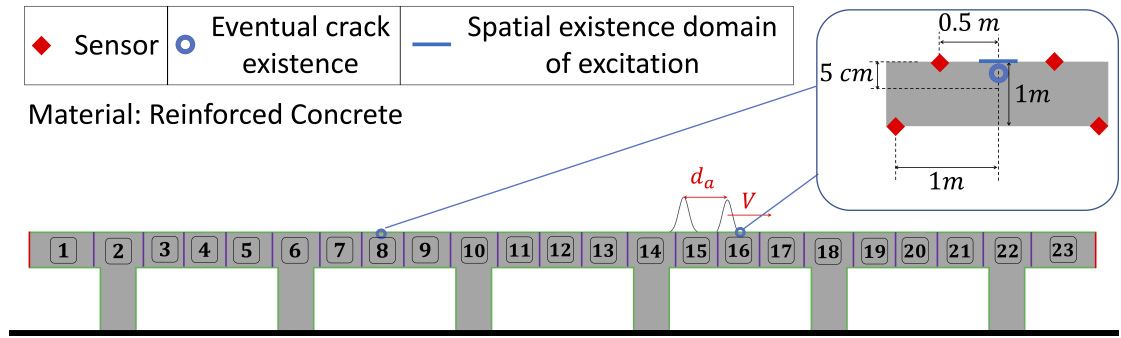}
	\caption{Note-to-scale representation of the global system for the bridge with different sensors location}\label{fig:sensors_further_bridge_SBC}	
\end{figure}

Sensors located further away from the potential cracks (centered at $(8\times L,H)$ and $(16\times L,H)$) are considered and the new sensors' positions are shown in figure \ref{fig:sensors_further_bridge_SBC} where the sensors belonging to the lower boundary of the decks have the same location as in the previous example, while the sensors located at the upper boundary are now $0.5 \ {\rm m}$ away from the potential cracks instead of the $0.2 \ {\rm m}$ chosen before. 

\begin{figure}[H]
	\centering
	\includegraphics[width=\textwidth]{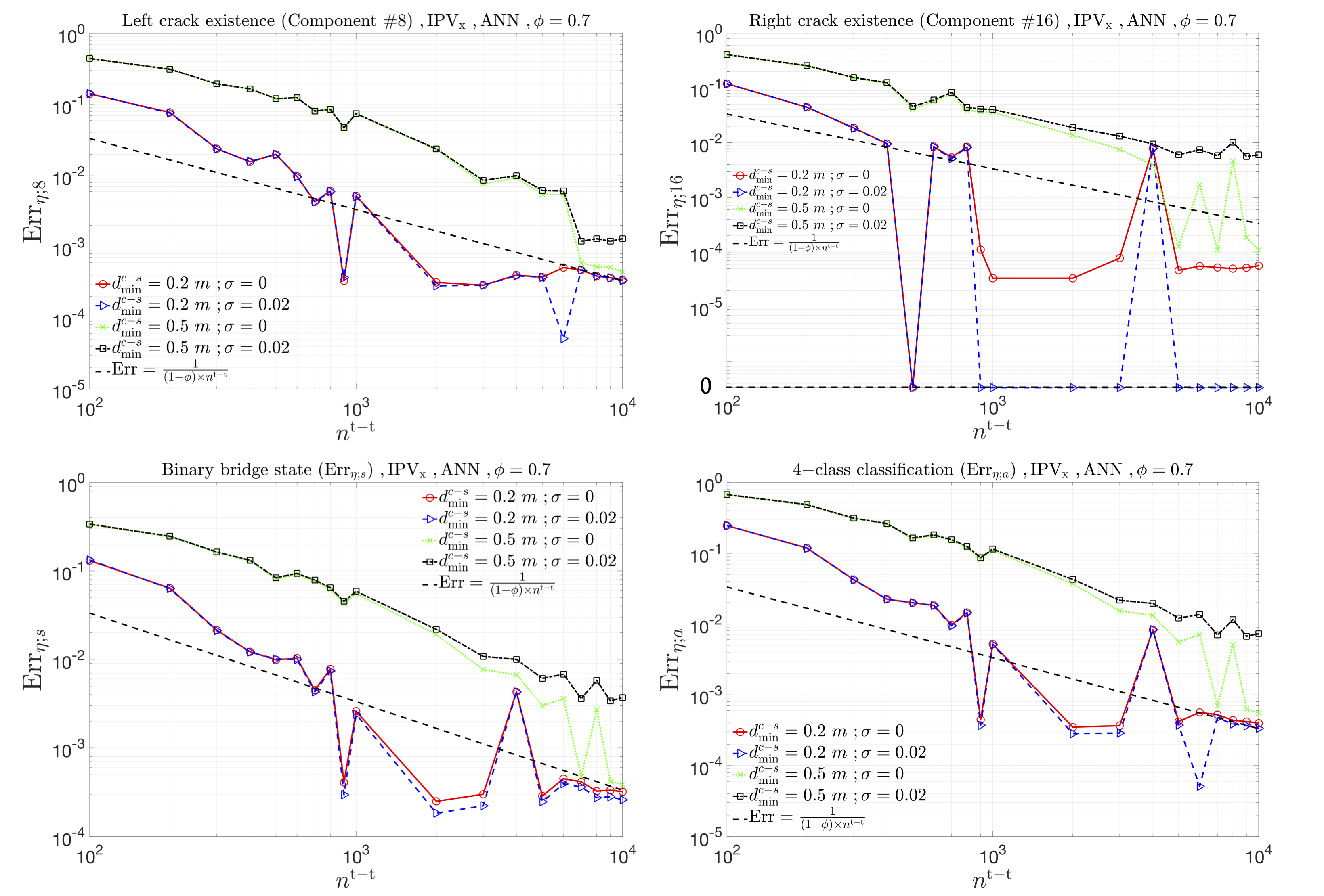} % 0.337
	\caption{Misclassification errors variation with $\ntt$ for different sensors locations using ${\rm IPV_x}$ features, ANN and with $\phi=0.7$: $d_{\rm min}^{c-s}$ refers to the minimum distance between the potential crack location and the closest sensor}\label{fig:bridge_classification_ANN_IPVx_new_comp_sensors_effect}	
\end{figure}

\begin{figure}[H]
	\centering
	\includegraphics[width=\textwidth]{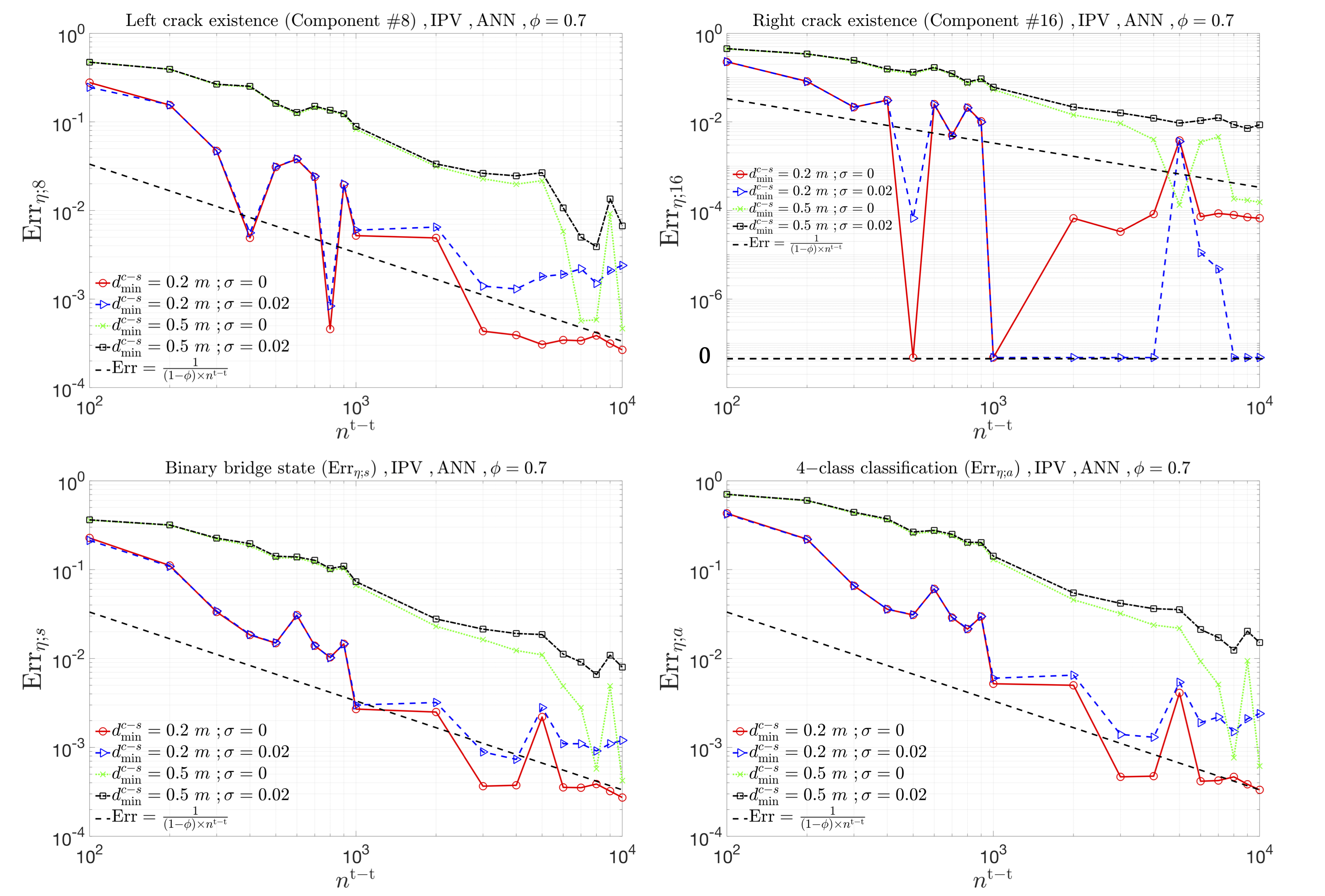} % 0.337
	\caption{Misclassification errors variation with $\ntt$ for different sensors locations using ${\rm IPV}$ features, ANN and with $\phi=0.7$: $d_{\rm min}^{c-s}$ refers to the minimum distance between the potential crack location and the closest sensor}\label{fig:bridge_classification_ANN_IPV_new_comp_sensors_effect}	
\end{figure}

Figure \ref{fig:bridge_classification_ANN_IPVx_new_comp_sensors_effect} shows the expected misclassifications obtained with the ${\rm IPV_x}$ features and ANN with $\sigma=0$ and $\sigma=0.02$ for the two different choices of sensors' locations. As expected, the misclassification error is higher for the new sensors' positions since they are located further away from the potential cracks. Nonetheless, for noiseless test dataset, the change in expected misclassification is not as significant as for test dataset with noise. Indeed, for $\sigma=0$ and as $\ntt$ gets close to $10^4$, the expected misclassification obtained for the new sensors' positions is on the same order as the one obtained for the sensors located closer to the cracks' locations. However, for test dataset altered with noise, the expected misclassification obtained for for sensors located further from the cracks is higher; and even for $\ntt$ close to $10^4$, the expected misclassifications are one order higher than those obtained for sensors located closer to the cracks. For instance, the expected misclassification obtained with the closest sensors located at $0.5 \ {\rm m}$ from the cracks is not less than one misclassified point out of the test dataset as $\ntt$ gets close to $10^4$, while it is the case when the closest sensors are located at $0.2 \ {\rm m}$ from the cracks. Nevertheless, the expected misclassifications always presents an almost monotonous decrease with the increase of the train-test dataset size $\ntt$. Hence, a longer distance between the sensors and the cracks along with a reasonable noise level requires larger train-test dataset, which emphasizes more the importance of pMOR and more specifically of the two-level PR-RBC approach in the context of large structure with localized excitations such that the sufficiently large dataset can be built in a reasonable computation time. The same conclusions can be drawn regarding the effect of the sensors' locations when considering the ${\rm IPV}$ features as shown in figure \ref{fig:bridge_classification_ANN_IPV_new_comp_sensors_effect}.

\begin{figure}[H]
	\centering
	\includegraphics[width=\textwidth]{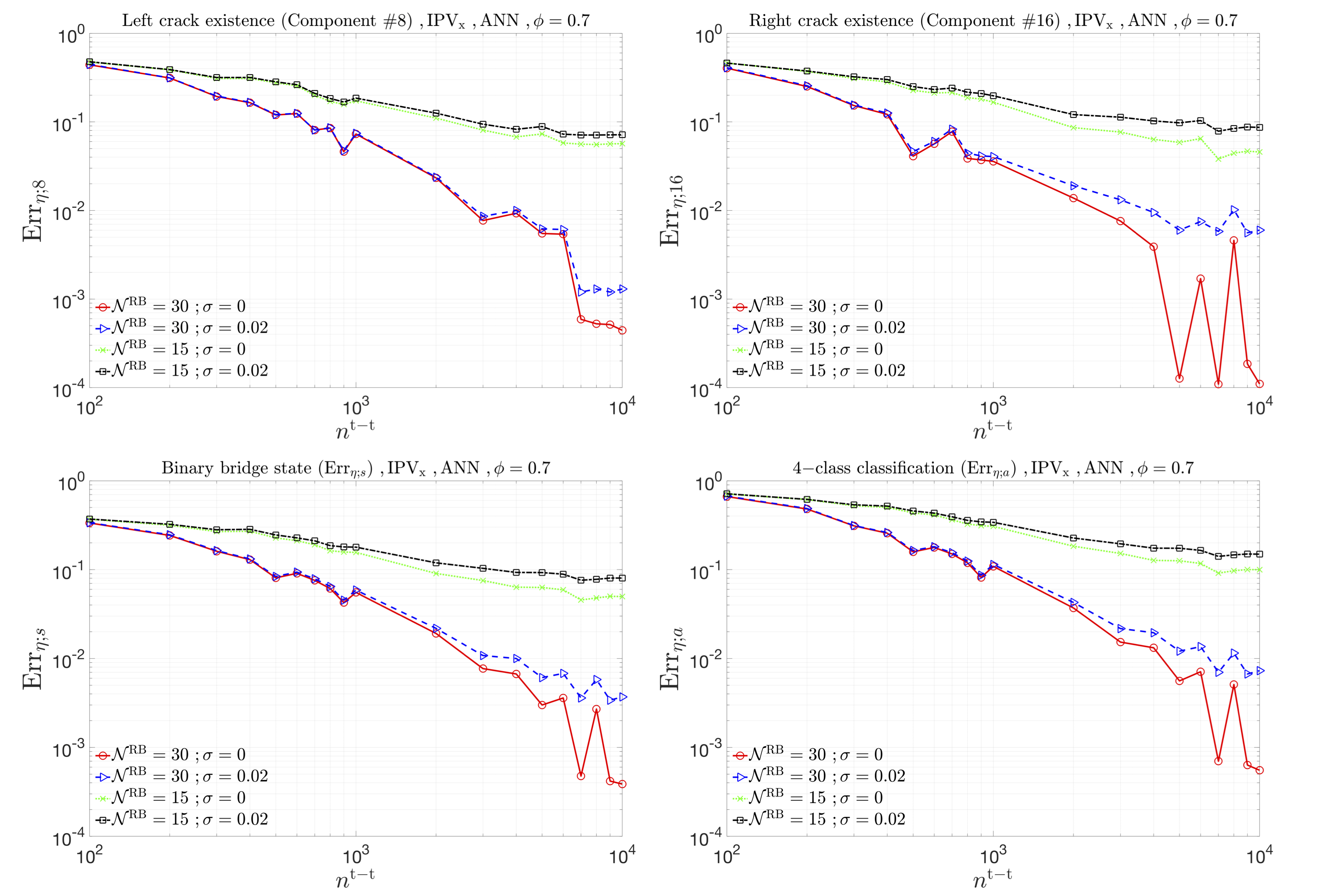}
	\caption{Effect of numerical solution accuracy on misclassification errors using ${\rm IPV_x}$ features, ANN and with $\phi=0.7$}\label{fig:bridge_classification_ANN_IPVx_new_comp_NRB_effect}	
\end{figure}

\begin{figure}[H]
	\centering
	\includegraphics[width=\textwidth]{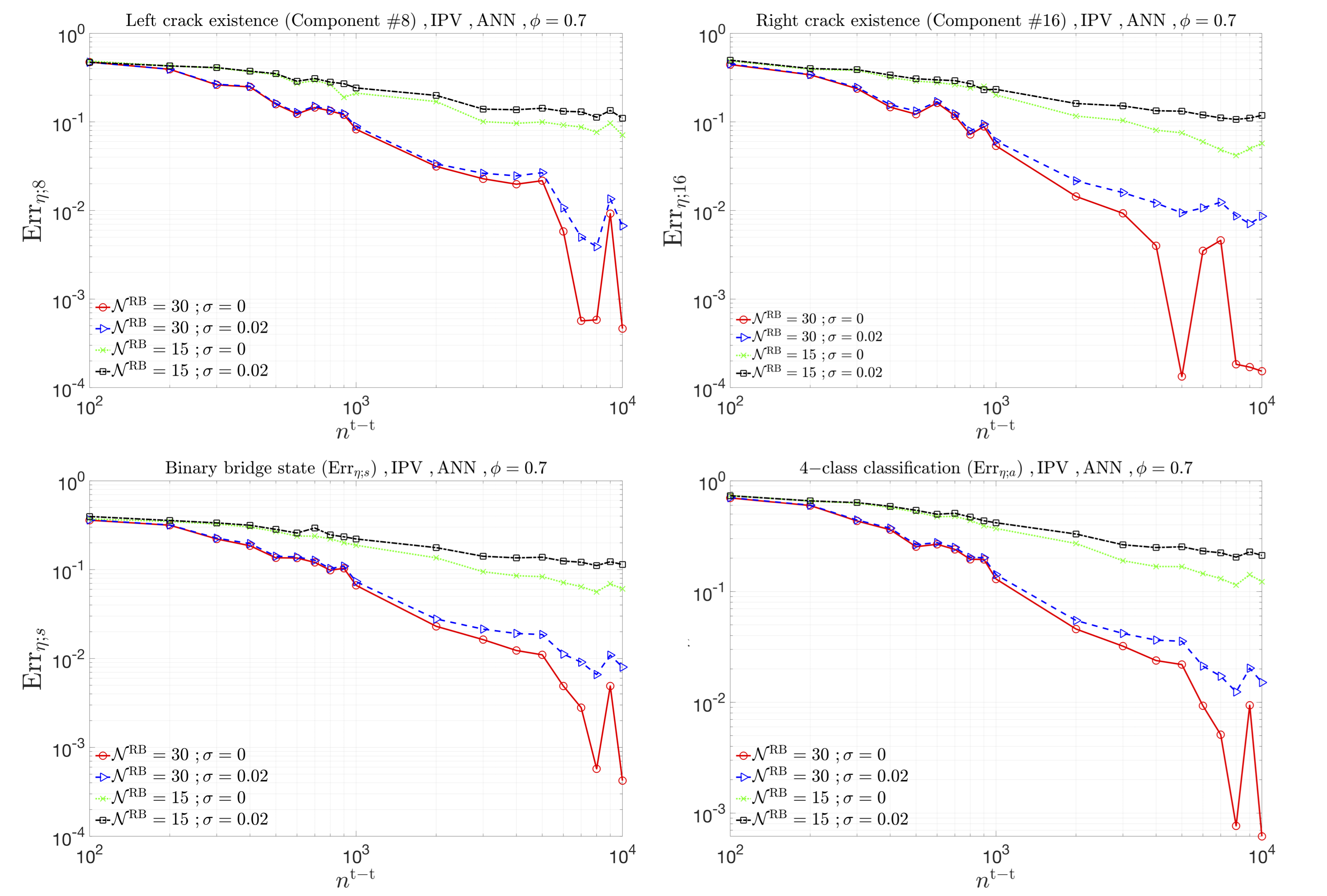}
	\caption{Effect of numerical solution accuracy on misclassification errors using ${\rm IPV}$ features, ANN, with $\phi=0.7$}\label{fig:bridge_classification_ANN_IPV_new_comp_NRB_effect}	
\end{figure}

Finally, in order to show its importance on the classification results, the accuracy of the two-level PR-RBC-based solution is intentionally reduced by, for instance, reducing the size $\calN^{\rm RB}$ of the final reduced space $X_{\rm RB}$ obtained at Level 2 reduction. The size needed to have a relative error in time of order of magnitude equal to $10^{-3}$ is around $30$ as shown in figure \ref{fig:err_FE_PR_RBC_by_max_time_bridge_PR_RBC}. If $\calN^{\rm RB}$ is limited to $15$, the relative error in time has an order of magnitude equal to $10^{-2}$. Figure \ref{fig:bridge_classification_ANN_IPVx_new_comp_NRB_effect} gives the misclassification results obtained with the ${\rm IPV_x}$ features and ANN, for $\sigma=0$ and $\sigma=0.02$, using the new sensors' locations (i.e. the closest sensor to any crack is $0.5 \ m$ away from it) and with the size of the final reduced space $\calN^{\rm RB}$ equal to $30$ and $15$. For any noise level considered, the misclassification accuracy suffers significantly when the two-level PR-RBC-based solution was intentionally made not accurate enough. The same conclusions can be drawn for the ${\rm IPV}$ features and ANN as shown in figure \ref{fig:bridge_classification_ANN_IPV_new_comp_NRB_effect}. %, and the same behavior was observed for any feature, any noise level and any classifier considered. Indeed, the misclassification accuracy suffers significantly when the two-level PR-RBC-based solution was made not accurate enough on purpose. Therefore, sufficiently accurate numerical solutions that can be computed within reasonable computation time to be able to build sufficiently large train-test synthetic dataset are required to have satisfactory classification results in the context of SHM for large structures with localized excitation, and the two-level PR-RBC method fulfills all these requirements.

%However, for $\calN^{\rm RB}=30$ and as stated above, as $\ntt$ gets close to $10^4$ the expected misclassification decreases such that there is statistically less than one misclassified point out of $3000$ ones, even for the bridge state classification task ($4$-class classification task). 

\begin{figure}[H]
	\centering
	\includegraphics[width=\textwidth]{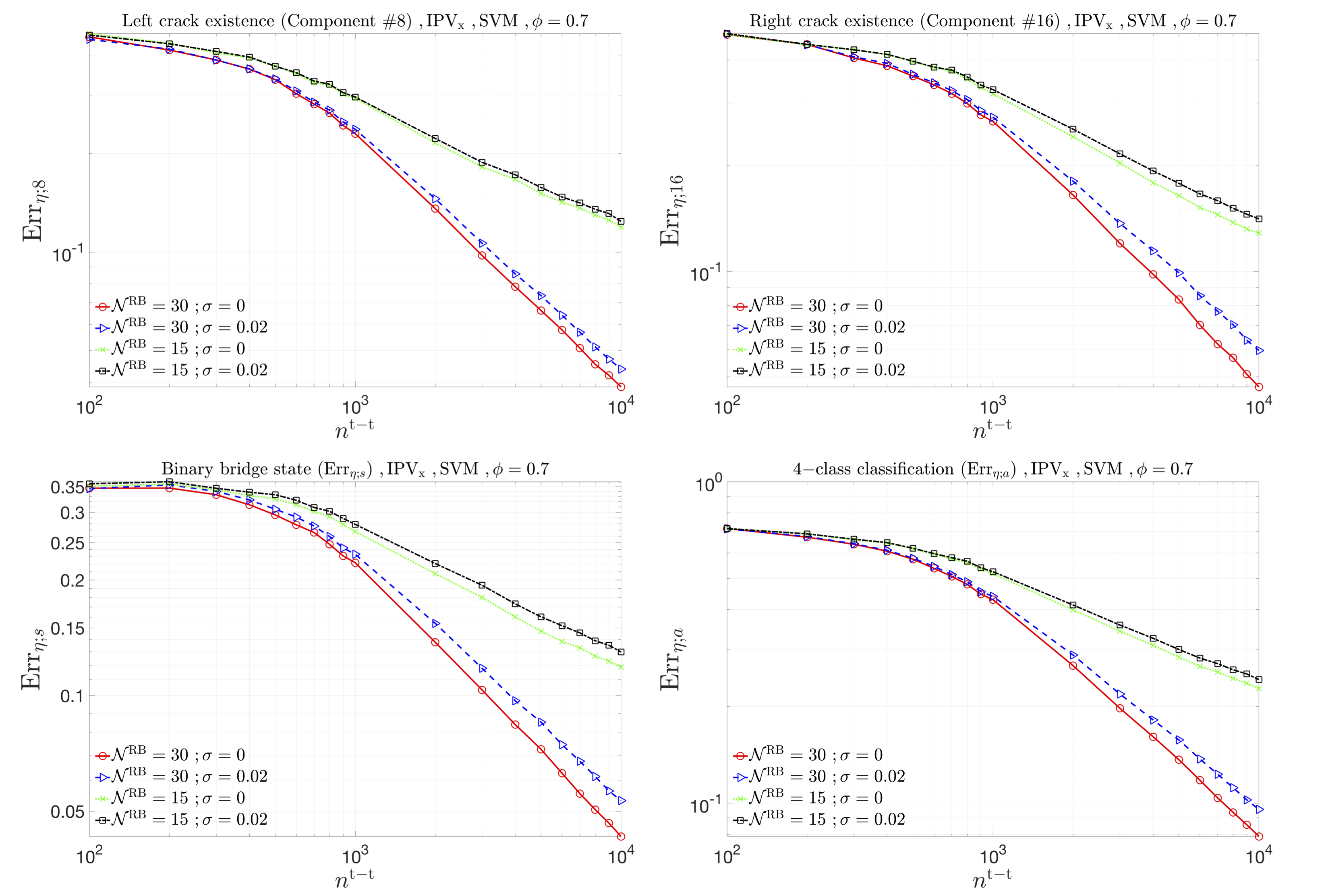}
	\caption{Effect of numerical solution accuracy on misclassification errors using ${\rm IPV_x}$ and SVM, with $\phi=0.7$}\label{fig:bridge_classification_SVM_IPVx_new_comp_NRB_effect}	
\end{figure}

The accuracy of the numerical solution is crucial for any feature, any noise level and any classifier considered. Indeed, figure \ref{fig:bridge_classification_SVM_IPVx_new_comp_NRB_effect} gives the misclassification results obtained with the ${\rm IPV_x}$ features and SVM, for $\sigma=0$ and $\sigma=0.02$, using the new sensors' locations (i.e. the closest sensor to the crack is $0.5 \ {\rm m}$ away from it) and with $\calN^{\rm RB}$ equal to $30$ and $15$. Again, for any noise level considered, the misclassification accuracy suffers significantly when the two-level PR-RBC-based solution was made not accurate enough on purpose. Therefore, sufficiently accurate numerical solutions that can be computed within reasonable computation time to be able to build sufficiently large train-test synthetic dataset are required to have satisfactory classification results in the context of SHM for large structures with localized excitation, and the two-level PR-RBC method fulfills all these requirements.

\section{Conclusions}\label{sec:Cl}
\subsection{Summary}

This work presents the development of a SBC approach for SHM of large deployed mechanical structures. The method takes advantage of the two-level PR-RBC method as a pMOR approach to approximate the hyperbolic PDE of time-domain elastodynamics for large domains with localized source terms and/or local parameters variation. The usefulness of such approach was demonstrated by considering a system with localized operational excitations and nuisance parameters. For such systems, frequency analyses fail to capture the structure response to such forces and hence the two-level PR-RBC approach developed for PDE of time-domain elastodynamics is of great usefulness. 
 %%% work from here %%%

For the simulation task, features based on the two-point correlation function were built and state-of-the-art machine learning algorithms were considered to perform a damage detection on the structure. The features considered are based on damage indices defined in previous works but which were never used within a context of SBC. Thus, the damage indices were already shown to be sensitive to the damage | thanks to the definition of the two-point correlation function | but the generalization of their former definition showed that the features used for the classification task in this work are relatively insensitive to nuisance parameters and measurement noise. Within the generalization, the reference sensor choice dilemma was solved since all possible combinations of the two-point correlation functions are considered and imbedded in the feature that will be sued for the classification. Moreover, instead of visual inspection of the damage indices, classifiers are exploited to automatically discriminate between the different damage cases. The two-level PR-RBC approach was used to construct a sufficiently large dataset with rich description of undamaged and damaged states. The system was modeled as a continuum governed by the elastodynamics PDE, providing a more faithful characterization than modeling the mechanical structure as an assembly of rigid beams. The two-level PR-RBC approach lets us considerably reduce the degrees of freedom of the numerical system and thus the computational cost related to the construction of the dataset. In addition, the model order reduction technique and the classification task inherently permit to account for the probabilistic nature of the nuisance parameters. Finally, thanks to the great flexibility in topology, the two-level PR-RBC approach can accommodate for more realistic damage instances such as crack existence and is thus particularly well-suited to SHM. Test classification errors below $0.1\%$ were reached for disjoint training set of size $7\times10^3$ and test set of size $3\times10^3$, thus showing the strong potential of the proposed approach in view of the application to real-life systems.

For the SBC task, a bridge example was considered in which the goal is to detect the existence of two potential cracks and such problem shows:
\begin{itemize}
\item the characterization of localized operational excitation in terms of nuisance parameters, 
\item the merits of time-domain-based correlation function features in the context of localized operational excitation, other nuisance variables and added noise, and
\item the importance of the two-level PR-RBC approach which considerably reduces the computational burden associated with the construction of synthetic training dataset. Thus, it allows the generation of such sufficiently large dataset with reasonable computation time, while guaranteeing the accuracy of the numerical solution. These two conditions need to be satisfied in order to obtain good classification results as shown by the significant decrease of the misclassification errors with the increase of the size of the training dataset.
\end{itemize} 

\subsection{Future Work}
As a future work, it would be interesting to develop computational techniques for automatic, or at least semi-automatic, feature identification which would include sensor placement and two-point correlation functions selection \cite{Qarib2015}. These methods promise not only more effective classifiers but also more effective deployment. For instance, they can reduce the number of sensors needed. To this end, recent advances in Robust Optimization \cite{Ben2009,Bertsimas2011} and Mixed Integer Optimization \cite{Bertsimas2005} can be useful. %These feature extraction techniques will be directly informed by the target risk and are in this sense ``optimal" with respect to our stated objectives. 

It would also be interesting to consider experimental data to use as a test for the classifiers that are built using synthetic dataset. Impact hammers represent an interesting option for SHM actuation in particular as regards installation and operation. Experimental enrichment of the training procedure could also be considered if some experimental data might be available at the offline stage of the classification task. %An example of impulsive SHM is Acoustic Pulse Reflectometry \cite{Amir2010,Sharp1997}: a microphone measures the reflections of an acoustic pulse in a pipeline to deduce damage ``scatterers" such as holes or obstructions; similar concepts may be applied more generally within the elastodynamics context. Experimental enrichment of the training procedure could also be considered if some experimental data might be available at the offline stage of the classification task. %We could also exploit techniques developed in Robust Optimization \cite{Ben2009,Bertsimas2011} to reduce the sensitivity of the classifier to data uncertainty.

Finally, in many applications it may not be feasible to provide an exhaustive experimental description of all possible states of damage. In other applications, a detailed characterization of damage may not in fact be interesting | but rather only a determination of ``undamaged" or ``damaged" (or, in the case of quality control, ``different") is required. From the learning or statistical perspective, this two-way classification corresponds to the problem of anomaly, or outlier, detection \cite{Chandola2009}. For the examples considered in this work, the classification results obtained for such binary task were presented. Within these two contexts, investigating unsupervised learning techniques would be of great interest.

%\backmatter

\textbf{Acknowledgement: }
%\section{Acknowledgments}
%\label{sec:ack}
This work was supported by the ONR Grant [N00014-17-1-2077] and by the ARO Grant [W911NF1910098]. We would like to thank Professor Anthony T. Patera, Dr. Tommaso Taddei and Professor Masayuki Yano for the helpful comments and software they provided us with.


\newpage
\bibliographystyle{elsarticle-num}
\nocite{*}% Show all bib entries - both cited and uncited; comment this line to view only cited bib entries;
\section*{\refname}
\bibliography{SBC_SHM_CO_ST.bib}

%%%%%\begin{thebibliography}{00}


%%%%%\bibitem{}

%%%%%\end{thebibliography}

\end{document}